\documentclass[12pt]{amsart}
\usepackage{diagrams}

\usepackage{amssymb}
\newcommand{\nc}{\newcommand}

\nc{\Ab}{{\mathtt {Ab}}}          
\nc{\Alg}{{\mathtt{Alg}}} 
\nc{\Algw}{{\mathtt{Algw}}} 

\nc{\art}{{\mathtt {art}}}
\renewcommand{\c}{\mathrm{c}}
\nc{\cdef}{{\CD}\mathit{ef}}
\nc{\cdel}{{\CD}\mathit{el}}
\nc{\Def}{\mathrm{Def}}
\nc{\Del}{\mathrm{Del}}
\nc{\cf}{\mathrm{cf}}
\nc{\cl}{\mathit{cl}}
\nc{\DO}{\mathit{do}}
\nc{\aff}{\mathrm{aff}}
\nc{\cff}{\mathrm{cff}}
\nc{\chom}{\mathrm{Hom}^{\bullet}}
\nc{\cphom}{\CH\mathit{om}^{\bullet}}
\nc{\eq}{\mathit{eq}}
\nc{\shom}{\mathrm{Hom}_{\bullet}}
\nc{\cpder}{{\CD\mathit{er}^{\bullet}}}
\nc{\der}{\mathrm{Der}}
\nc{\Zar}{\mathrm{Zar}}

\nc{\cocart}{\mathrm{cocart}}
\nc{\Cone}{\mathrm{cone}}
\nc{\Coind}{\mathrm{Coind}}

\nc{\Coker}{\mathrm{Coker}}
\nc{\cosk}{\mathrm{cosk}}
\nc{\diag}{\mathrm{diag}}
\nc{\Dop}{\Delta^{\op}}
\nc{\f}{\mathrm{f}}
\nc{\fl}{{f\!\ell}}
\nc{\Hom}{\operatorname{Hom}}
\nc{\invlim}{\underset{\leftarrow}{\lim}\,}
\nc{\iso}{\mathrm{iso}}
\nc{\isom}{\overset{\sim}{\lra}}
\nc{\Ind}{\mathrm{Ind}}
\nc{\Ker}{\mathrm{Ker}}
\nc{\Mod}{{\mathtt{Mod}}}
\nc{\mult}{\mathrm{mult.}}       
\nc{\Nerve}{\CN}
\nc{\Nil}{{\mathtt {Nil}}}
\nc{\Norm}{\mathrm{Norm}}
\nc{\pt}{\mathrm{pt}}

\nc{\op}{{\mathrm{op}}}

\nc{\phom}{\CH\!o\!m}
\nc{\pre}{\begin{picture}(0.2,0.4)
          \put(0,-0.2){$\widehat{\ }$}
          \end{picture}}

\nc{\Vectgr}{{\mathtt {Vectgr}}} 
\nc{\sCat}{\mathtt{sCat}}
\nc{\sGrp}{\mathtt{sGrp}}
\nc{\she}{\sim}
\nc{\sHo}{\mathtt{sHo}}
\nc{\sHoalg}{\mathtt{sHoalg}}
\nc{\simpl}{\Dop\Ens}
\nc{\sk}{\mathrm{sk}}
\nc{\SEQ}{\CS\!\mathit{eq}}
\nc{\WEQ}{\CW\!\mathit{eq}}
\nc{\TKS}{T^{\mathrm{KS}}}
\renewcommand{\vert}{\mathrm{vert}}
\nc{\hor}{\mathrm{hor}}

\nc{\invtensor}{\underset{\leftarrow}{\otimes}}
\nc{\rlarrows}{\begin{picture}(1,0.4)
                \put(0,-0.1){\makebox(1,0.2){$\leftarrow$}}
                \put(0,0.1){\makebox(1,0.2){$\to$}}
              \end{picture}}
\nc{\rra}{\begin{picture}(1,0.4)
                \put(0,-0.1){\makebox(1,0.2){$\lra$}}
                \put(0,0.1){\makebox(1,0.2){$\lra$}}
              \end{picture}}
\nc{\Left}{\mathbf L}  
\nc{\Right}{\mathbf R} 
\nc{\gr}{\operatorname{gr}}
\nc{\Ho}{\operatorname{Ho}}
\nc{\alt}{\operatorname{alt}}
\nc{\Sym}{\operatorname{Sym}}
\nc{\sym}{\operatorname{sym}}
\nc{\id}{\operatorname{id}}
\nc{\im}{\operatorname{Im}}
\nc{\Col}{\operatorname{Col}}
\nc{\ter}{\operatorname{ter}}
\nc{\intl}{\operatorname{int}}
\nc{\val}{\operatorname{val}}

\nc{\TN}{{\cal N}}
\nc{\Nor}{\operatorname{N}}
\nc{\Tor}{\operatorname{Tor}}
\nc{\res}{\operatorname{res}}
\nc{\Stab}{\operatorname{Stab}}

\nc{\End}{\operatorname{End}}
\nc{\holim}{\operatorname{holim}}
\nc{\dirlim}{\underset{\rightarrow}{\lim}\,}
\nc{\com}{\operatorname{co}}
\nc{\Tot}{\operatorname{Tot}}
\nc{\Th}{\operatorname{Th}}
\nc{\Cech}{\check{C}}
\nc{\Spec}{\operatorname{Spec}}
\nc{\Spf}{\operatorname{Spf}}
\nc{\MC}{\operatorname{MC}}
\nc{\U}{\operatorname{U}}
\nc{\Diff}{{\cal D}\mbox{\em iff}}
\nc{\Mor} {{\cal M}or}
\nc{\Ob}{\operatorname{Ob}}
\nc{\cone}{\widehat}
\nc{\Coder}{\operatorname{Coder}}
\nc{\pr}{\operatorname{pr}}

\nc{\CHo}{{\cal{H}\mbox{\it{o}}}}
\nc{\Modf}{{\mathtt{modf}}}       
\nc{\Modg}{{\mathtt{modg}}}       
\nc{\Hoalg}{{\mathtt {Hoalg}}} 
\nc{\Valg}{{\mathtt {Viral}}} 
\nc{\Algf}{{\mathtt {Algf}}} 
\nc{\Algg}{{\mathtt {Algg}}}

\nc{\dgc}{{\mathtt{dgc}}}
\nc{\dgca}{{\mathtt{dgca}}}
\nc{\dgcu}{{\mathtt{dgcu}}}
\nc{\dgcuf}{{\mathtt{dgcuf}}}
\nc{\dgcf}{{\mathtt{dgcf}}}
\nc{\dgcg}{{\mathtt{dgcg}}}
\nc{\dgcc}{{\mathtt{dgccc}}}

\nc{\dgl}{{\mathtt{dglie}}}
\nc{\dgla}{{\mathtt{dgla}}}
\nc{\dglf}{{\mathtt{dglf}}}
\nc{\dglg}{{\mathtt{dglg}}}
\nc{\dga}{{\mathtt{dga}}}
\nc{\dgar}{{\mathtt {dgart}^{\leq 0}}}
\nc{\Coll}{{\mathtt{Coll}}}
\nc{\Kan}{{\mathtt {Kan}}}

\nc{\Grp}{{\mathtt {Grp}}}
\nc{\Groups}{{\mathtt {Groups}}}

\nc{\Cat}{{\mathtt {Cat}}}
\nc{\Ens}{{\mathtt {Ens}}}
\nc{\Op}{{\mathtt{Op}}}

\nc{\Lie}{{\mathtt{LIE}}}
\nc{\Com}{{\mathtt{COM}}}
\nc{\Ass}{{\mathtt{ASS}}}

\nc{\pa}{\partial}

\nc{\G}{{\boldmath{G}}}

\nc{\cal}{\mathcal} 

\nc{\CA}{\cal A}
\nc{\CBB}{\cal B}
\nc{\CB}{\cal B}
\nc{\CC}{\cal C}
\nc{\CD}{\cal D}
\nc{\CE}{\cal E}
\nc{\CF}{\cal F}
\nc{\CG}{\cal G}
\nc{\CH}{\cal H}
\nc{\CI}{\cal I}
\nc{\CJ}{\cal J}
\nc{\CK}{\cal K}
\nc{\CL}{\cal L}
\nc{\CM}{\cal M}
\nc{\CN}{\cal N}
\nc{\CO}{\cal O}
\nc{\CP}{\cal P}
\nc{\CQ}{\cal Q}
\nc{\CR}{\cal R}
\nc{\CS}{\cal S}
\nc{\CT}{\cal T}
\nc{\CU}{\cal U}
\nc{\CV}{\cal V}
\nc{\CW}{\cal W}
\nc{\CZ}{\cal Z}

\nc{\fa}{\mathfrak a}
\nc{\fg}{\mathfrak g}
\nc{\fk}{\mathfrak k}
\nc{\fh}{\mathfrak h}
\nc{\fm}{\mathfrak m}
\nc{\fn}{\mathfrak n}
\nc{\fA}{\mathfrak A}
\nc{\fC}{\mathfrak C}
\nc{\fF}{\mathfrak F}
\nc{\fI}{\mathfrak I}
\nc{\fP}{\mathfrak P}
\nc{\fS}{\mathfrak S}

\nc{\nen}{\newenvironment}
\nc{\ol}{\overline}
\nc{\ul}{\underline}
\nc{\lra}{\longrightarrow}
\nc{\lla}{\longleftarrow}
\nc{\Lra}{\Longrightarrow}
\nc{\Lla}{\Longleftarrow}
\nc{\Llra}{\Longleftrightarrow}
\nc{\hra}{\hookrightarrow}

\nc{\notebox}[1]{\noindent\fbox{\parbox{12.5cm}{\sf #1}}\\[8pt]}
\nc{\Thm}[1]{Theorem~\ref{#1}}
\nc{\Prop}[1]{Proposition~\ref{#1}}
\nc{\Lem}[1]{Lemma~\ref{#1}}
\nc{\Cor}[1]{Corollary~\ref{#1}}
\nc{\Conj}[1]{Conjecture~\ref{#1}}
\nc{\Claim}[1]{Claim~\ref{#1}}
\nc{\Defn}[1]{Definition~\ref{#1}}
\nc{\Exa}[1]{Example~\ref{#1}}
\nc{\Rem}[1]{Remark~\ref{#1}}
\nc{\Note}[1]{Note~\ref{#1}}

\nen{thm}[1]{\label{#1}{\bf Theorem.\ } \em}{}
\nen{prop}[1]{\label{#1}{\bf Proposition.\ } \em}{}
\nen{lem}[1]{\label{#1}{\bf Lemma.\ } \em}{}
\nen{klem}[1]{\label{#1}{\bf Key Lemma.\ } \em}{}
\nen{cor}[1]{\label{#1}{\bf Corollary.\ } \em}{}
\nen{conj}[1]{\label{#1}{\bf Conjecture.\ } \em}{}
\nen{claim}[1]{\label{#1}{\bf Claim.\ } \em}{}

\nen{defn}[1]{\label{#1}{\bf Definition.\ } }{}
\nen{exa}[1]{\label{#1}{\bf Example.\ } }{}

\nen{rem}[1]{\label{#1}{\em Remark.\ } }{}
\nen{note}[1]{\label{#1}{\em Note.\ } }{}
\nen{exer}[1]{\label{#1}{\em Exercise.\ } }{}
\nen{pf}{\begin{proof}}{\end{proof}}

\begin{document}

\title[Deformations of sheaves]
{Deformations of sheaves of algebras}
\author{Vladimir Hinich}
\address{Department of Mathematics, University of Haifa,
Mount Carmel, Haifa 31905,  Israel}
\email{hinich@math.haifa.ac.il}
\begin{abstract}A construction of the tangent dg Lie algebra
of a sheaf of operad algebras on a site is presented. The
requirements on the site are very mild; the requirements on the
algebra are more substantial. A few applications including
the description of deformations of a scheme and equivariant 
deformations are considered.
The construction is based upon a model structure on the category
of presheaves which should be of an independent interest.
\end{abstract}
\maketitle
%
\setcounter{section}{-1}
\section{Introduction}

%
%

\subsection{}
In this paper we study formal deformations of sheaves of algebras. The most 
obvious (and very important) example is that of deformations of a scheme
$X$ over a field $k$ of characteristic zero. In two different cases, the first
when $X$ is smooth, and the second when $X$ is affine, the description is 
well-known. In both cases there is a differential graded (dg) Lie algebra
$T_X$ over $k$ such that formal deformations of $X$ over the artinian local 
base $(R,\fm)$ are described by the Maurer-Cartan elements of 
$\fm\otimes T_X$, modulo a gauge equivalence.

It is well-understood now that formal deformations over a field of
characteristic zero are governed by a differential graded Lie algebra.
One of possible explanations of this phenomenon was
suggested in~\cite{uc}: we expect deformation problems to have formal moduli
(which is expected to be a ``commutative'' formal dg scheme). Then the 
representing dg Lie algebra corresponds to the formal moduli by Koszul
(or bar-cobar) duality. Thus, the existence of dg Lie algebra governing
deformations is equivalent to the representability (in ``higher'', dg sense)
of the deformation problem.

However, in the two cases mentioned above ($X$ smooth and $X$ affine) the 
governing dg Lie algebra $T_X$ appears in seemingly different ways. 
This can be shortly described as follows.

\subsubsection{$X$ is smooth}
\label{smooth}

Affine smooth scheme $X$ has no formal deformations. Its trivial
deformation $U_R$ with an artinian local base $(R,\fm)$ 
admits the automorphism group $\exp(\fm\otimes T_X)$ which is nothing
but the value at $R$ of the formal group corresponding to 
$T_X=\Gamma(X,\CT_X)$. Descent theorem of~\cite{ddg} asserts in this 
situation that for a general smooth scheme $X$ the dg Lie algebra $T_X$
governing the deformations of $X$ can be calculated by the formula
\begin{equation*}
\label{T.smooth}
T_X=\Right\Gamma(X,\CT_X).
\end{equation*}

\subsubsection{$X$ is affine}
\label{affine}

Let $X=\Spec(A)$ for a commutative $k$-algebra $A$. It is convenient
to consider $A$ as a dg commutative $k$-algebra concentrated at degree zero.
Then the deformation theory of dg algebras ~\cite{dha} suggests the following
recipe of calculation of $T_X$. Let $P\to A$ be a cofibrant (some call
it free or semi-free) resolution of $A$ in the model category of 
commutative dg $k$-algebras. Deformations of $A$ and of $P$ are equivalent; 
deformations of $P$ appear as perturbations of the differential which 
are described by the Maurer-Cartan elements of the Lie algebra of 
derivations of $P$. Thus, one has
\begin{equation*}
\label{T.affine}
T_X=\der(P,P).
\end{equation*}

\subsubsection{} We wish to describe in a similar way deformations of a 
sheaf of algebras. The first problem seems to be the lack of cofibrant
resolutions for sheaves of algebras. This turns out to have a very pleasant
solution: the category of complexes of presheaves admits a model category
structure describing the homotopy theory of complexes of sheaves.

A similar model category structure exists for sheaves of operad algebras
in characteristic zero. This allows us to define deformation functor
and to construct the corresponding dg Lie algebra in a way similar to
the one described in~\ref{affine}. The construction is local, so, as a 
result, we obtain a presheaf of dg Lie algebras. We use the construction 
mentioned in~\ref{smooth} to get a global dg Lie algebra. 

\subsection{Sheaves vs presheaves}

Let $X$ be a site and $F:X^{\op}\to \CA$ be a presheaf on $X$ with values in
a category $\CA$ having a notion of weak equivalence (for instance,
complexes, simplicial sets, categories or polycategories). The notion of 
sheaf is not very appropriate here: we know this well, for instance, in the 
case $\CA=\Cat$. This was probably the reason  Jardine~\cite{j} suggested
a model category structure on the category of simplicial presheaves. The idea
was to extend the notion of weak equivalence so that a presheaf will
be weakly equivalent to its sheafification. Then the localization of the
category of simplicial sheaves with respect to the weak equivalences 
can be described as the homotopy category of the category of presheaves.

We adopt a similar point of view. We need a model category structure on 
presheaves of algebras which would allow us to construct 
``semi-free resolutions''. This model category structure is based upon a model
category structure on the category $C(X^{\pre}_k)$ of complexes of 
presheaves of $k$-modules which is described by the following result.

\subsubsection{}
\begin{thm}{intro:cmc-com}
Let $X$ be a site, $k$ a ring and let $C(X^{\pre}_k)$ denote the
category of complexes of presheaves of $k$-modules on $X$.

1. The category $C(X^{\pre}_k)$ admits a model category structure so that
\begin{itemize}
\item{ weak equivalences are maps $f:M\to N$ inducing a quasi-isomorphisms
$f^a:M^a\to N^a$ of sheafifications.}
\item{ cofibrations are generated by maps $f:M\to M\langle x; 
dx=z\in M(U)\rangle$ coresponding to adding a section to kill a cycle $z$
over an object $U\in X$.}
\end{itemize}

2. A map $f:M\to N$ in the above model category structure is a fibration
iff $f(U):M(U)\to N(U)$ is surjective for any $U\in X$ and for any 
hypercover $\epsilon: V_{\bullet}\to U$ of $U\in X$ the corresponding 
commutative diagram
\begin{center}
\begin{diagram}
M(U)& \rTo & \Cech(V_{\bullet},M)       \\
\dTo&      &  \dTo                      \\
N(U)& \rTo & \Cech(V_{\bullet},N)
\end{diagram} 
\end{center}
is homotopy cartesian.
\end{thm}

We remind the notion of hypercover in~\ref{HC}. Cech complex
$\Cech(V_{\bullet},M)$ of $M$ with respect to a hypercover $V_{\bullet}$
is defined as the total complex corresponding to the cosimplicial
complex $n\mapsto M(V_n)$, see~\ref{cech-compl}.

Notice that we do not require the existence of limits in the site $X$.
This is important for us since we want to be able to apply this to 
the category of affine open subsets of a scheme which usually does not admit a 
final object. 

\subsubsection{}The proof of~\Thm{intro:cmc-com} is given in~\ref{pf-cmc-com}.
It is based on an explicit description of generating acyclic cofibrations.

Recently (see~\ref{ack}) we learned that the model structure described
above (at least part 1 of \Thm{intro:cmc-com}) is known to  specialists,
see, for instance, \cite{toen}, Appendix C. We decided, however, to present 
our proof since it is direct, general, and gives an explicit description of 
fibrations which we need in any case. A similar model category structure was 
used in~\cite{hs} for $n$-stacks. 

We also present in Appendix~\ref{App:sim} a version of \Thm{intro:cmc-com} for 
simplicial presheaves. This model category structure on simplicial 
presheaves has the same weak equivalences as Jardine's \cite{j}
but the cofibrations are generated by gluing cells and fibrations 
have a similar description using hypercovers. 

Under some mild restrictions on the $X$ (any hypercover can be refined
by a split hypercover) our model structure coincides with a one 
recently defined in~\cite{dhi}, Theorem 1.3.

\subsection{Higher deformation functor}

Classical formal deformation functors can be usually described as follows.
Let $\art(k)$ be the category of artinian local $k$-algebras
$(R,\fm)$ with residue field $k$. Let $\CC(R)$ for each $R\in\art(k)$
denote  a groupoid of ``objects over $R$'', so that for each map $f:R\to S$ 
in $\art(k)$ a base change functor $\alpha^*:\CC(R)\to \CC(S)$ is defined.

Then the groupoid of formal deformations of an object $A\in\CC(k)$ over 
$R\in\art(k)$ is defined as the fiber of $\pi^*:\CC(R)\to\CC(k)$ at $A$.

In higher deformation theory one extends the category $\art(k)$ allowing
artinian algebras which are not necessarily concentrated at degree zero.
In this paper we work with the category $\dgar(k)$ of non-positively 
graded differential artinian algebras, see~\cite{uc} for the explanation.
One cannot expect that the deformation functor extended in this way 
has values in the category of groupoids: one should expect 
a higher version of groupoid appearing here. 
We use simplicial groupoids (or simplicial weak groupoids
which is the same from the homotopical point of view) as a higher
version of groupoid. 

Our higher formal deformation functors can be described as follows.

Let $\dgar(k)$ be the category of non-positively graded commutative 
artinian local dg algebras over $k$ with residue field $k$.
Let $\CC(R)$ for each $R\in\dgar(k)$ denote a category of ``objects over $R$''.
We suppose there is a subcategory $\CW(R)$ of weak equivalences in $\CC(R)$.
Let $\widehat{\CW}(R)$ be the full Dwyer-Kan (hammock) localization of $\CW(R)$
(see \ref{weak-grp} for the details). This is a simplicial groupoid and we
 define the simplicial groupoid of formal deformations of $A\in\CC(k)$  
as the homotopy fiber of the map 
\begin{equation*}
\widehat{\CW}(R)\to\widehat{\CW}(k)
\end{equation*}
at $A$.
  
For the description of deformations of sheaves of algebras we take
$\CC(R)$ to be the category of sheaves of $R$-algebras flat over $R$.
$\CW(R)$ is the subcategory of quasi-isomorphisms of sheaves of algebras. 

\subsection{Main result}

Let $X$ be a site and let $k$ be a field of characteristic zero.
Let  $\CO$ be an operad in the category of complexes of sheaves of 
$k$-modules on $X$ and let $A$ be a sheaf of $\CO$-algebras.

Our main result, \Thm{result}, presents (under some restrictions on $X$, 
$\CO$ and $A$)  the dg Lie algebra governing formal deformations of $A$. 
The construction goes as follows. According 
to~\Thm{cmc-alg}, the category of presheaves of $\CO$-algebras admits
a model structure generalizing the one defined in~\ref{intro:cmc-com}.
Let $P$ be a fibrant cofibrant $\CO$-algebra weakly equivalent to $A$.
Then the presheaf of derivations of $P$, $\CT_A:=\cpder(P,P)$, is
a fibrant presheaf of dg Lie algebras on $X$. The dg Lie algebra
governing deformations of $A$ can be expressed then as $\holim\CT_A$.

Here are the assumptions for which the result is proven.

\subsubsection{Assumptions on $X$} 
\begin{itemize}
\item{The topos $X^{\she}$  admits enough points.}
\item{The final presheaf in $X^{\pre}$  admits a finite hypercover. }
\end{itemize}
The second condition is of course fulfilled for sites
admitting a final object. However, in our main application $X$ is the site
of affine open subschemes of a scheme. In this case the condition
is fulfilled for quasi-compact separated or finite dimensional schemes.

\subsubsection{Assumptions on $\CO$} 

Complexes $\CO(n)$ are  non-positively graded.
This condition does not seem to be really restrictive. Operads one
encounters are usually obtained by a tensor product of a sheaf of rings 
(e.g, the structure sheaf) with a constant operad.

\subsubsection{Restrictions on $A$} 

For each $U\in X$ the cohomology $H^i(U,A)$ is supposed to vanish for  $i>0$. 
This condition means the following. Choose a fibrant
resolution $A\to A'$ of $A$. Then for each $U\in X$ the complex $A'(U)$
has no positive cohomology. This is the most serious assumption. It is not
in general fulfilled even in the case $A$ is a sheaf of algebras. In fact,
in this case the cohomology $H^i(U,A)$ has its usual meaning (cohomology 
of the sheaf $A|_U$) and it does not vanish in general for any $U\in X$.

The situation is, however, slightly better then one could think. The reason
is that once we are given a sheaf of algebras $A$ in a topos $X^{\she}$,
we have a freedom in the choice of $X$. If $X^{\she}$ admits a generating
family of sheaves $U$ satisfying $H^i(U,A)=0$, we can choose $X$ to be the 
site generated by this family. For instance, if $A$ is a quasi-coherent sheaf
on a scheme, one chooses $X$ to be the category of affine open subschemes of 
the scheme, with the Zariski topology. 

\subsection{Applications}

Direct application of \Thm{result} gives
the following result (see \ref{schemes} and \ref{qcoh}).

{\em
Let $X$ be a scheme over a field of characteristic zero. Suppose $X$
admits a finite dimensional hypercover by affine open subschemes. Then 
the functor of formal deformations of $X$ (or, more generally, of a
quasicoherent operad algebra on $X$) is represented by a dg Lie algebra. 
}

The dg Lie algebra representing deformations of a quasicoherent sheaf
of algebras, is usually difficult to determine. Its cohomology, however,
can be easily identified with the Hochschild cohomology (for associative
algebras), see \ref{obstruction.theory}. In a very special case
of associative deformations of the structure sheaf of a smooth scheme,
the tangent Lie algebra identifies with the (shifted and truncated)
complex of Hochschild cochains given by polydifferential operators.

The last application we present in this paper is to the description of 
equivariant deformations. Let $A$ be a sheaf of algebras on a site $X$
satisfying the conditions of~\Thm{result}, and let $T$ be 
the dg Lie algebra governing the deformations of $A$. Suppose now that
a formal group $G$ acts on $X$ and on $A$ in a compatible way. Then
$G$ acts in a natural way on $T$ and the $G$-equivariant deformations 
of $A$ are governed by a dg Lie algebra $\Right\Gamma^G(T)$ whose
$i$-th cohomology is $H^i(G,T)$. This is proven 
in~\ref{equivariant.deformations}.

\subsection{Structure of the sections}

In Section~\ref{sect:models} we prove \Thm{intro:cmc-com} describing the
model category structure on the category of complexes of presheaves.
We describe the functors $\Right\cphom$ and $\Right\Gamma$ using this
model category structure. In Section~\ref{sect:modop} we present a model
structure for the category of presheaves of operad algebras. In
Section~\ref{sect:result} we describe the deformation functor for
sheaves of operad algebras on a site. Here the main \Thm{result}
is proven. In Section~\ref{sect:examples}
we present two examples: deformations of a scheme (or of a quasi-coherent
algebra on a scheme) and equivariant deformations of a sheaf of algebras
with respect to a discrete group.

In Appendix~\ref{App:scat} we present a necessary information about simplicial
categories, Dwyer-Kan localization and its different presentations
for a model category. In Appendix~\ref{App:sim} which is not used in the main
body of the paper, we present a model category structure on the category
of simplicial presheaves and provide a description of fibrations
similar to that of~\Thm{intro:cmc-com}.

\subsection{Notation}
In this paper $\mathbb{N}$ denotes the set of non-negative integers,
$\Delta$ the category of sets $[n]=\{0,\ldots,n\},\ n\in\mathbb{N}$
and of non-decreasing maps, and $\Ens$ denotes the category of sets.
The category of simplicial sets is denoted $\Delta^{\op}\Ens$.

As well, we denote $\Ab$ the category of abelian groups, $\Cat$ the category
of small categories, $\Grp$ the subcategory of groupoids.

If $\CA$ is an abelian category, $C(\CA)$ is the category of complexes over 
$\CA$; if $\CA$ is a tensor category, $\Op(\CA)$ is the category of
operads in $\CA$. The notation $\Mod$ and $\Alg$ for the categories 
of modules and algebras is obvious. $\Sigma_n$ is the symmetric group.

\subsection{Relation to other works} This work extends the approach of
\cite{dha} to the sheaves of algebras. Both 
\cite{dha} and the present work are based
on an idea (which goes back to Halperin-Stasheff~\cite{hast},
Schlessinger-Stashef~\cite{schst}, Felix~\cite{fel}) that deformations
of an algebra can be described by perturbation of the differential in its
free resolution. Since ~\cite{hast,schst,fel} a better understanding
of the notion of deformation has been achieved, due to Drinfeld and 
Deligne, so that the language of obstructions is being
substituted with the dg Lie algebra formulation of deformation theory. 

One has to mention Illusie~\cite{illu} and Laudal~\cite{laud} who 
constructed the obstruction theory for deformation of schemes, and 
Gerstenhaber-Schack ~\cite{gesch} who studied obstruction theory for 
presheaves of algebras. Obstruction theory for deformations of sheaves 
of associative algebras was studied in~\cite{gait} and~\cite{lunts}.

\subsection{Acknowledgements}
A part of this work was made during my visits 
at MPIM and at IHES. I am grateful to these institutions for stimulating
atmosphere and excellent working conditions. During the conference on 
polycategories at Nice (November, 2001) I knew that a part of the 
results on model category structures described here is known to 
specialists. I am very grateful to the organizers of the conference 
A.~Hirshowitz, C.~Simpson, B.~Toen for the invitation.

\section{Models for sheaves}
\label{sect:models}

%
%
%

Let $X$ be a site, $X^{\pre}$ and $X^\she$ be the categories of
presheaves (resp., sheaves) on $X$. If $k$ is a commutative ring,
$X^{\pre}_k$ (resp., $X^\she_k$) denotes the category or 
preseaves (resp., sheaves) of $k$-modules on $X$.

The categories $X^{\pre}_k$ and $X^\she_k$ are tensor (=symmetric 
monoidal) categories. The sheafification functor $M\mapsto M^a$ is exact and
preserves the tensor product, see~\cite{sga4}, IV.12.10.

In this section we provide a model (=closed model category) 
structure for the category $C(X^{\pre})$ of complexes of presheaves on $X$. 
This structure ``remembers'' the topology of $X$ in a way that that the 
model category  $C(X^{\pre})$ becomes a powerful tool in doing 
homological algebra of {\em sheaves} on $X$.

In the next section we will describe a similar model structure
on the category of presheaves of algebras over a dg operad on $X$.

\subsection{Coarse topology}

\subsubsection{}
\begin{thm}{coarse-mod}
The category $C(X^{\pre}_k)$ of presheaves of $k$-modules admits a model
structure with weak equivalences defined as pointwise quasi-isomorphisms 
and fibrations as pointwise surjections.
\end{thm}

Since the category $C(k)$ of complexes of $k$-modules is cofibrantly 
generated (see ~\cite{dhk}, 7.4, for the definition of cofibrantly 
generated model categories and~\cite{haha} for the model structure 
on $C(k)$), the result follows from the general observation 
of~\cite{dhk}, 9.6. The CMC structure described above is also 
cofibrantly generated.

Let us recall the description of a generating collection of cofibrations.
Let $U\in X$. We will identify $U$ with the presheaf represented by $U$.
The generating cofibration given by a pair $(U\in X,n\in\mathbb{Z})$ is 
defined as
\begin{equation*}
i:\, kx\cdot U\to ky\cdot U\oplus kz\cdot U
\end{equation*}
where $n=|x|=|y|,\ n-1=|z|,\ i(x)=y,\ dx=dy=0,\ dz=y.$

The generating acyclic cofibration is defined for each 
$(U\in X,n\in\mathbb{Z})$ as

\begin{equation}
\label{coarse-gac}
j:\, 0\to kx\cdot U\oplus ky\cdot U
\end{equation}
where $n=|x|,\ n-1=|y|,\ dx=0,\ dy=x.$

The model category structure defined in \ref{coarse-mod} knows nothing
about the topology of $X$. It corresponds to the coarse topology on $X$.

In the general case the notion of hypercover is of a great importance.

\subsection{Hypercovers}
\label{HC}

Let us recall a few standard notions connected to hypercovers.
The context presented here is slightly more general than that of~\cite{sga4},
Expos\'e V.

Let $X$ be a site. Notice that we do not require that fiber products and finite
products exist in $X$.

\subsubsection{}
An object $K\in X^{\pre}$ is semi-representable if it is isomorphic to 
a coproduct of representable presheaves.

A simplicial presheaf $K_{\bullet}$ is called {\em hypercover} of $X$ if

\begin{itemize}
\item[(HC0)] {For each $i\geq 0$ $K_i$ is semi-representable.}
\item[(HC1)] { For each $n\geq 0$ the canonical map 
$$ K_{n+1}\to (\cosk_n(K))_{n+1}$$
is a cover of presheaves (i.e. its sheafification is surjective).}
\item[(HC2)] {The canonical map of presheaves $K_0\to *$ is a cover.}
\end{itemize}

Let $L\in X^{\pre}$. 
A simplicial presheaf $K_{\bullet}$ endowed with an augmentation
$\epsilon:K_0\to L$ is called a hypercover of $L$ if it defines a
hypercover of the site $X/L$.

\subsubsection{}
Let $K_{\bullet}$ be a simplicial presheaf on $X$. We denote by
$C_*(K,k)$, or simply $C_*(K)$, the complex of normalized chains of $K$.
This is an object of $C(X^{\pre}_k)$. If $L$ is a presheaf of sets considered
as a disrete simplicial presheaf, $C_*(L)=kL$ is the free presheaf of vector
spaces generated by $L$.

The following lemma is of crucial importance for us.

\begin{lem}{sga4lem}{\em (cf.~\cite{sga4}, V.7.3.2(3)).}
Let $\epsilon:K_{\bullet}\to L$ be a hypercover. Then the induced
map $C_*(K_{\bullet})\to kL$ induces a quasi-isomorphism of sheafifications. 
\end{lem}

\subsection{General case}
\label{cmc-gen-com}
Now we will define another CMC structure on $C(X^{\pre}_k)$ which 
``remembers'' about the topology on $X$.

\subsubsection{}
\begin{thm}{gen-com}
1. The category $C(X^{\pre}_k)$ of presheaves of $k$-modules admits a CMC 
structure with cofibrations as in~\ref{coarse-mod} and weak equivalences 
defined as maps $f:M\to N$ such that the sheafification $f^a$ is a  
quasi-isomorphism of complexes of sheaves. 

2. A map $f:M\to N$ is a fibration iff $f(U):M(U)\to N(U)$
is surjective for each $U\in X$ and 
for any 
hypercover $\epsilon: V_{\bullet}\to U$ of $U\in X$ the corresponding 
commutative diagram
\begin{center}
\begin{equation}
\label{diag:com-fib}
\begin{diagram}
M(U)& \rTo & \Cech(V_{\bullet},M)       \\
\dTo&      &  \dTo                      \\
N(U)& \rTo & \Cech(V_{\bullet},N)
\end{diagram} 
\end{equation}
\end{center}
is homotopy cartesian.
\end{thm}

\subsubsection{}
\label{pf-cmc-com}
The definition of Cech complex corresponding to a hypercover is given 
in~\ref{cech-compl}. The first part of~\Thm{gen-com} is proven 
in~\ref{gener-ac}--\ref{end1:gen-com}. It is based on an expicit 
description of the collection of generating acyclic cofibration. 
The second part of the theorem is proven in~\ref{end2:gen-com}.

\subsubsection{}
\label{gener-ac}
The generating set of acyclic cofibrations is numbered by pairs 
$$(\epsilon:V_{\bullet}\to U,n)$$
where $\epsilon$ is a hypercover and $n$ an integer. An acyclic cofibration
$j:K\to L$ corresponding to a pair $(\epsilon,n)$ as above is defined as 
follows.

The presheaf $K$ is shifted by $-n$ cone of the mapping
$$ C_*(\epsilon):C_*(V_{\bullet})\to U.$$

This means that $K^n=k\cdot U$ and 
$K^{n-i-1}=k\cdot V_i/\sum_j k\cdot s_jV_{i-1}$ for $i\geq 0$,
where $s_j:V_{i-1}\to V_i,\quad j=0,\ldots,i-1$ are the degenerations.

The presheaf $L$ is defined as the cone of $\id_K$, with the obvious 
canonical embedding $j:K\to L$.

We will write sometimes $K_{\epsilon,n}$ and $L_{\epsilon,n}$ for 
the complexes $K$ and $L$ corresponding to a pair 
$(\epsilon:V_{\bullet}\to U,n)$.

\subsubsection{}
\begin{lem}{ac-is-we}
The map $j:K\to L$ of presheaves constructed above is a weak equivalence.
\end{lem}
\begin{proof}
The map $j$ is obviously injective. The sheafification of its cokernel
is contractible by \cite{sga4}, V.7.3.2(3).
\end{proof}

\subsubsection{}
\begin{note}{a-is-con}
The sheafification of $K$ is also contractible by 
\cite{sga4}, V.7.3.2(3). 
\end{note}

\subsubsection{}
\begin{lem}{af}
Let a map $f: M\to N$ of presheaves satisfy the right lifting property
with respect to all generating acyclic cofibrations and let 
$f^a$ be a (pointwise) quasi-isomorphism. Then $f$ is pointwise surjective
quasi-isomorphism.
\end{lem}
\begin{proof}
The map $f$ is pointwise surjective by~\ref{coarse-mod}. Therefore, 
we have to prove that for any $U\in X$ the map  $f(U)$\, is a 
quasi-isomorphism. We can put $N=0$ without loss of generality. 
Let $a\in M(U)^n$ be a cycle. 
We have to prove $a$ is a boundary. For this we will construct a 
hypercover $\epsilon: V_{\bullet}\to U$ and a map
\begin{equation*}
f: K_{\epsilon,n}\to M
\end{equation*}
whose restriction to the $n$-th component is the map $a:U\to M^n$.
Then by the right lifting property $f$ lifts to a map $L_{\epsilon,n}\to M$
which proves $a$ is a boundary.

We proceed by induction: since $M^a$ has zero cohomology
and $da=0$, there exists a cover $\epsilon: V_0\to U$ such that 
$\epsilon^*(a)$ is a boundary. Fix $a_0\in M(V_0)$ such that 
$da_0=\epsilon^*(a)$. Suppose, by the induction hypothesis, there exist
sections $a_i\in M(V_i),\quad i=0,\ldots,n$ such that
$a_i$ vanishes on the degeneracies $M(s_jV_{i-1})$ and
$da_i=\sum\ (-1)^j d_j^*(a_{i-1})$. Then one can choose a cover
$V_{n+1}\to\cosk_n(V_{\bullet})_{n+1}$ for which the cycle
$\sum\ (-1)^id_i^*(a_n)$ becomes a boundary. 
\end{proof}

\subsubsection{}
\label{end1:gen-com}
Let $J$ be the collection of generating acyclic cofibrations. We define
fibrations as the maps of presheaves satisfying the RLP with respect to the
elements of $J$. We denote by $\ol{J}$ the collection of maps which can be 
obtained as a countable direct composition of pushouts of coproducts of maps 
in $J$. We call the maps from $\ol{J}$ standard acyclic cofibrations. They are
cofibrations and weak equivalences by~\Lem{ac-is-we}.
 
 According to \Lem{af}, fibrations which are weak equivalences
are precisely acyclic fibrations in the sense of~\ref{coarse-mod} (i.e., 
pointwise surjective quasi-isomorphisms). 

Let $f:A\to B$ be a map of presheaves. The existence of decomposition $f=pi$
where $p$ is an acyclic fibration and $i$ is a cofibration follows from 
 \Thm{coarse-mod}. A small object argument, see~\cite{dhk}, II.7.3,
implies the existence of a decomposition $f=qj$ where $j\in \ol{J}$ and
$q$ is a fibration. Suppose now $f$ is a cofibration and a weak equivalence
and choose a decomposition $f=qj$ as above. The map $q$ is therefore a 
fibration and a weak equivalence, therefore, by~\Lem{af}, $q$ is an 
acyclic fibration in the sense of~\ref{coarse-mod}. Therefore, $q$ satisfies
the RLP with respect to $f$. This implies that $f$ is a retract of $j$.

We proved that any acyclic cofibration is a retract of a standard acyclic 
cofibration which yields the first part of the theorem. The second part 
of the theorem is explained in~\ref{end2:gen-com} below.

\subsubsection{}
\label{cech-compl}
Let $V_{\bullet}$ be a simplicial presheaf and let $M\in C(X^{\pre}_k)$.

The collection 
\begin{equation}
\label{cos-cech}
n\mapsto\Hom(V_n,M)
\end{equation}
is a cosimplicial object in $C(X^{\pre}_k)$. \v{C}ech complex of $M$, 
$\Cech(V_{\bullet},M)$, is defined as the normalization
of~(\ref{cos-cech}), so that
\begin{equation}
\label{cech}
\Cech^n(V_{\bullet},M)=\{f\in\prod_{p+q=n}\Hom(V_p,M^q)|f\text{ vanishes
on the degenerate simplices }\}.
\end{equation} 
In particular, any hypercover $\epsilon:V_{\bullet}\to U$ gives rise
to a map 
$$ M(U)\to\Cech(V_{\bullet},M).$$

Note that one has an obvious isomorphism
\begin{equation}
\label{cech-vs-chain}
\Cech(V_{\bullet},M)=\chom(C_*(V_{\bullet}),M).
\end{equation} 

\subsubsection{Proof of \Thm{gen-com}(2)}
\label{end2:gen-com}
Any fibration is pointwise surjective since the maps~(\ref{coarse-gac})
are acyclic cofibrations. From now on we suppose that $f$ is pointwise
surjective. Put $K=\Ker(f)$. 
The condition~(\ref{diag:com-fib}) is now equivalent to the claim that
the map 
$$ K(U)\to \Cech(V_{\bullet},K)$$
is a weak equivalence. This is equivalent to the requirement that
the complex $ \chom(K_{\epsilon,0},K)$ has trivial cohomology or, equivalently,
that the map
$$ \chom(K_{\epsilon,0},M)\to\chom(K_{\epsilon,0},N)$$
is a quasi-isomorphism.
 This, in turn, can be interpreted as  
the right lifting property of $f$ with respect to the generating
acyclic cofibrations $K_{\epsilon,n}\to L_{\epsilon,n}$.
%
%
\qed

\subsubsection{}
\begin{note}{dependence-on-site}
The model category structure in $C(X^{\pre}_k)$ described in the
theorem, depends essentially on the site $X$ (i.e., on the generating
family of the topos $X^{\she}$). For instance, a quasi-coherent sheaf 
is not fibrant in the Zariski site of a scheme. However, it is fibrant
when considered as a presheaf on the site of {\em affine} open subschemes of 
the scheme.
\end{note}

The following observation will be useful in the sequel.

\subsubsection{}
\begin{lem}{restriction}
Let $X$ be a site and let $U\in X$. Let $j:X/U\to X$ be the natural embedding.
The restriction functor $j^*: C(X^{\pre}_k)\to C((X/U)^{\pre}_k)$
preserves fibrations and weak equivalences. If $X$ admits finite products,
$j^*$ preserves cofibrations.
\end{lem}
\begin{proof}
Preservation of weak equivalences is immediate. Preservation of fibrations
follows immediately from~\Thm{gen-com}(2). To prove that $j^*$
preserves cofibrations we will check that the right adjoint functor
$j_*: C((X/U)^{\pre}_k)\to C(X^{\pre}_k)$ preserves acyclic fibrations.
One has
\begin{equation*}
j_*(M)(W)=M(U\times W)
\end{equation*}
for $M\in C((X/U)^{\pre}_k)$;
acyclic fibrations are just pointwise surjective quasi-isomorphisms.
This implies the lemma.
\end{proof}

\subsection{Cohomology and $\Right\Gamma$}
\label{cohomology}

\subsubsection{}
Let $M\in C(X^{\pre}_k)$ and let $\widetilde{M}$ be a fibrant 
resolution of $M$. We define the $i$-th cohomology presheaf of $M$, 
$\CH^i(M)$, to be the presheaf
\begin{equation*}
\label{i-th.coh}
\CH^i(M)(U)=H^i(\widetilde{M}(U)).
\end{equation*}
The following lemma shows the result does not depend of the choice of the 
resolution.

\subsubsection{}
\begin{lem}{}Let $f:M\to N$ be a weak equivalence of fibrant objects
in $C(X^{\pre}_k)$. Then for each $U\in X$ and $i\in\mathbb{Z}$
the map $H^i(M(U))\to H^i(N(U))$ is bijective.
\end{lem}
\begin{proof}
Any acyclic fibration induces a pointwise quasi-isomorphism, so we can suppose
that $f$ is an acyclic cofibration. Since $M$ is fibrant, $f$ is split by a 
weak equivalence $g:N\to M$. Therefore, the map $H^i(f(U))$ is injective.
This reasoning can be also applied to the map $g$ instead of $f$, to conclude
that $H^i(g(U))$ is injective as well.
\end{proof}

\subsubsection{Global sections} 
\label{RGamma} 
Suppose first that $X$ admits a final object
$*\in X$. Then the derived global sections $\Right\Gamma(M)$ can be defined as
$\widetilde{M}(*)$ where $\widetilde{M}$ is a fibrant resolution of $M$.

If $X$ does not admit a final object, one defines the derived global sections 
functor as follows. Choose a hypercover  $V_{\bullet}$ of $X$. We define
\begin{equation*}
  \label{eq:RG}
  \Right\Gamma(M)=\Cech(V_{\bullet},\widetilde{M}),
\end{equation*}
see~\ref{cech-compl}.
By formula~(\ref{cech-vs-chain}) and~\ref{a-is-con},
the result does not depend on the choice of hypercover $V_{\bullet}$ 
and of the fibrant resolution $\widetilde{M}$.

\subsection{Presheaves of modules}
Let $\CO$ be a presheaf of dg associative $k$-algebras on $X$. We denote \,
$\Mod(\CO,X^{\pre})$\, the category of (presheaves of dg) $\CO$-modules.

\subsubsection{}
\begin{thm}{O-mod}The category \, $\Mod(\CO,X^{\pre})$\, admits a model
structure for which a map $f:M\to N$ of presheaves of $\CO$-modules 
is a weak equivalence (resp., a fibration) iff it is a weak 
equivalence (resp., a fibration) of presheaves of $k$-modules.
\end{thm}

The proof of the theorem is easily deduced from the following lemma.

\subsubsection{}
\begin{lem}{h1-o-mod}
Let $j:K\to L$ be a generating acyclic cofibration corresponding
to a pair $(\epsilon,n)$ as in~\ref{gener-ac},
$M$ be a  $\CO$-module and $f:K\to M$ be a map of complexes of presheaves. 
Then the induced map
\begin{equation*}
\label{gac-o-mod}
\textstyle
M\to M\coprod^{\CO\otimes K}\CO\otimes L
\end{equation*} 
is a weak equivalence.
\end{lem}
\begin{proof}
The map in question being injective, it is enough to study the cokernel
which is isomorphic (up to a shift) to $\CO\otimes\Coker(j)$.
Its sheafification is isomorphic to $\CO^a\otimes\Coker(j)^a$.
The complex $\Coker(j)^a$ is acyclic by \Lem{sga4lem}.
Since its components are flat, the tensor product is acyclic as well.
\end{proof}

The following lemma shows that weakly equivalent associative algebras 
give rise to equivalent derived categories of modules.

\subsubsection{}Let now an algebra homomorphism $f:\CO\to \CO'$ be given.
One defines in a standard way a pair of adjoint functors
\begin{equation*}
  \label{modtomod}
f^*:\Mod(\CO,X^{\pre})\rlarrows\Mod(\CO',X^{\pre}):\ f_*
\end{equation*}
which induces a pair of derived functors

\begin{equation}
  \label{Dmodtomod}
\Right f^*:\Mod(\CO,X^{\pre})\rlarrows\Mod(\CO',X^{\pre}):\ f_*=\Left f_*.
\end{equation}

\begin{lem}{eq:Dmodtomod}
Let $f:\CO\to \CO'$ be a weak equivalence of presheaves of associative algebras.
Then the adjoint pair~(\ref{Dmodtomod}) establishes an equivalence of the
derived categories of modules.
\end{lem}
\begin{proof}
One has to check that of $M$ is a cofibrant $\CO$-module then the natural map
$$M\to f_*(f^*(M))$$
is a weak equivalence.
The claim immediately reduces to the case $M=\CO\otimes k\cdot U$ for 
$U\in X$. Then $f_*(f^*(M))=\CO'\otimes k\cdot U$ and the claim is obvious.
\end{proof}

\subsection{Inner $\cphom$}
The model category $\Mod(\CO,X^{\pre})$ admits an extra structure similar to 
that of simplicial model category of Quillen ~\cite{ha}.

\subsubsection{}
\begin{defn}{ihom}Let $M,\, N\in\Mod(\CO,X^{\pre})$. The inner Hom
$$ \cphom_{\CO}(M,N)\in C(X^{\pre}_k)$$
assigns to each $U\in X$ the complex of $k$-modules defined as
\begin{equation*}
\label{cphom}
\cphom_{\CO}(M,N)(U)=\underset{V\to V'\to U}{\invlim}
\chom_{\CO(V')}(M(V'),N(V)).
\end{equation*}
\end{defn}   
Here $\chom_{\CO(V')}$ is the usual inner Hom in the category of complexes 
of $\CO(V')$-modules. We will write as well $\chom_{\CO}(M,N)$ for the
complex of the global sections of $\cphom_{\CO}(M,N)$.

\subsubsection{}
\begin{lem}{IMC}Let $\alpha:M\to M'$ be a cofibration 
and $\beta:N\to N'$ be a fibration in $\Mod(\CO)$. Then the natural
map
\begin{equation*}
\label{map:IMC}
\cphom(M',N)\to\cphom(M,N)\times_{\cphom(M,N')}\cphom(M',N')
\end{equation*}
is a fibration in $C(X^{\pre}_k)$. It is a weak equivalence if $\alpha$
or $\beta$ is a weak equivalence.
\end{lem}

Standard adjoint associativity isomorphism
\begin{equation*}
\chom_{\CO}(X\otimes_kY,Z)\isom\chom_k(X,\cphom_{\CO}(Y,Z))
\end{equation*}
reduces the claim to the following.

\subsubsection{}
\begin{lem}{IMC2} (Here $\otimes=\otimes_k$). 
Let $\alpha:A\to A'$ be a cofibration in $C(X^{\pre}_k)$ 
and let $\beta:M\to M'$ be a cofibration in $\Mod(\CO,X^{\pre})$.
Then the induced map
\begin{equation}
\label{map:IMC2}
\textstyle
A\otimes M'\coprod^{A\otimes M}A'\otimes M\to A'\otimes M'
\end{equation}   
is a cofibration. It is an acyclic cofibration if $\alpha$ or $\beta$ is.
\end{lem}
\begin{proof}
For the first claim it is enough to check the case $\alpha$ and $\beta$
are generating cofibrations. This is a very easy calculation. For the
second claim note that the cokernel of~(\ref{map:IMC2}) is isomorphic
to $\Coker(\alpha)\otimes\Coker(\beta)$. Since the sheafification commutes
with the tensor product and $\Coker(\alpha)$ is flat, the result is immediate.

\end{proof}

\subsection{Comparing to sheaves}

Let $\CO$ be a dg algebra in $X^{\she}_k$.

Recall that one can define the functor $\Right\cphom_{\CO}$
on the category $\Mod(\CO,X^{\she})$ using Spaltenstein's notion
of $K$-injective complex of sheaves~\cite{sp}. A complex 
$I\in\Mod(\CO,X^{\she})$ is called $K$-injective if it satisfies the RLP
with respect to injective quasi-isomorphisms of sheaves.

Thus, one defines $\Right\cphom_{\CO}(M,N)$ as $\cphom_{\CO}(M,I)$
where $N\to I$ is a $K$-injective resolution.

\subsubsection{}
\begin{lem}{innerhoms} Here $\CO,\ M,\ N$ are as above.
Let $M'\to M$ be a cofibrant resolution of $M$ in $\Mod(\CO,X^{\pre})$
and $N\to N'$ be a fibrant resolution of $N$  in $\Mod(\CO,X^{\pre})$.
Then $\cphom_{\CO}(M',N')$ and $\Right\cphom_{\CO}(M,N)$ are equivalent.
\end{lem}
\begin{proof}
Let $N\to I$ be a $K$-injective resolution of $N$.
Any $K$-injective complex is fibrant as a complex of presheaves.
Therefore, $\cphom_{\CO}(M',N')$ and $\cphom_{\CO}(M',I)$ are weakly
equivalent fibrant complexes of presheaves. On the other hand,
$\cphom_{\CO}(M',I)=\cphom_{\CO}((M')^a,I)\isom\cphom_{\CO}(M,I)$
since $I$ is $K$-injective.
\end{proof}

\subsubsection{}
\begin{rem}{}
Subsection \ref{cohomology} and \Lem{innerhoms} show that the standard 
homological  algebra of sheaves can be rephrased in the language of 
complexes of presheaves endowed with the model structure defined 
in~\Thm{O-mod}. This has an advantage over the standard approach with
sheaves since the analog of \Thm{O-mod} takes place for presheaves
of algebras as well.
\end{rem}

\subsubsection{}
We have to mention the following consequence of~\ref{innerhoms}.

Let $M$ be a cofibrant and $N$ be a fibrant object in $\Mod(\CO,X^{\pre})$.
Let $U\in X$ and $j:X/U\to X$ be the natural embedding. According 
to \Lem{restriction} the functor $j^*$ preserves fibrations and weak 
equivalences, but does not necessarily preserve cofibrations. However,
the following takes place.

\begin{prop}{restriction:OK}
The complex $\cphom_{\CO|_U}(j^*(M),j^*(N))$ ``calculates the $\Right\cphom$''.
More precisely, if $M'\to j^*(M)$ is a cofibrant resolution, then
the induced map
$$ \cphom_{\CO|_U}(j^*(M),j^*(N))\to\cphom_{\CO|_U}(M',j^*(N))$$
is a weak equivalence. 
\end{prop}

\begin{proof}
Let $N\to I$ be a $K$-injective resolution of $N$. Then 
$j^*N\to j^*I$ is a $K$-injective resolution of $j^*N$.

This implies that all morphisms in the composition below
are equivalences.

\begin{multline}
\cphom_{\CO|_U}(j^*(M),j^*(N))=
j^*\cphom_{\CO}(M,N)\to
j^*\cphom_{\CO}(M,I)=\\
\cphom_{\CO|_U}(j^*(M),j^*(I))\to  
\cphom_{\CO|_U}(M',j^*(I))\leftarrow
\cphom_{\CO|_U}(M',j^*(N)).
\end{multline}

\end{proof}

\section{(Pre)sheaves of operad algebras}
\label{sect:modop}

%
%

%
%
%
%

In this section we describe a model  structure on the
category of presheaves of algebras over a $\Sigma$-split operad,
in the case when the corresponding topos has enough points. 
The structure is based on the model structure on the category
of complexes of presheaves described in \Thm{gen-com}.

We also discuss the category of modules over an operad algebra, derivations and
modules of differentials.

\subsection{Homotopical amenability}
Mimicing ~\cite{vir}, Definition 2.2.1, we define homotopically amenable 
presheaves of operads. 

\subsubsection{Notation}
Let $X$ be a site and let $k$ be a commutative ring. A presheaf
of operads on $X$ is, by definition, an operad in the tensor category
$C(X^{\pre}_k)$ of complexes of presheaves of $k$-modules.

If $\CO\in\Op(C(X^{\pre}_k))$ is such a presheaf, we denote by 
$\Alg(\CO,X^{\pre})$ the category of (presheaves of) $\CO$-algebras.

\subsubsection{}
\begin{defn}{hamo}
An operad $\CO\in\Op(C(X^{\pre}_k))$ is called {\em homotopically
amenable} if for each $A\in\Alg(\CO,X^{\pre})$ and for each
generating acyclic cofibration $j:K\to L$ in $C(X^{\pre}_k)$ 
with $K=K_{\epsilon,n}$ and $L=L_{\epsilon,n}$
where $\epsilon$ is a hypercover of $U\in X$ and $n\in\mathbb{Z}$,
see~\ref{gener-ac}, 
the map $\alpha$ defined by the cocartesian diagram
\begin{equation}
\label{defn:hamo}
  \begin{diagram}
     F(\CO,K) & \rTo^{j} & F(\CO,L) \\
      \dTo    &               &  \dTo    \\
       A      & \rTo^{\alpha}  &  A' 
  \end{diagram}
\end{equation}
is a weak equivalence.
\end{defn}

Here $F(\CO,K)$ denotes the free $\CO$-algebra generated by $K$. 

\subsubsection{}
The following result is standard.

\begin{thm}{cmc-alg:hamo}
Let $\CO$ be a homotopically amenable operad in $C(X^{\pre}_k)$. 
Then the category
$\Alg(\CO,X^{\pre})$ of $\CO$-algebras admits a model structure with weak
equivalences and fibrations defined as for $C(X^{\pre}_k)$ 
in~\ref{gen-com}.
\end{thm}

\subsection{$\Sigma$-split operads}
In the case $X^{\pre}=\Ens$  $\Sigma$-split operads defined in~\cite{haha},
 4.2, are homotopically amenable. Recall that a $\Sigma$-splitting
of an operad $\CO$ is defined as a collection of $\Sigma_n$-equivariant
splittings of the canonical maps
$$ \CO(n)\otimes k\Sigma_n\to \CO(n),\quad o\otimes\sigma\mapsto o\sigma$$
satisfying some extra compatibility properties, see~\cite{haha}, 4.2.4
for the precise definition.

The definition of $\Sigma$-split operad makes sense in any tensor category.
Thus, we can speak about $\Sigma$-split presheaves of operads on $X$.

It is worthwhile to mention two big classes of $\Sigma$-split operads.

\begin{itemize}
 \item{ If $k\supseteq\mathbb{Q}$ then all operads in $C(X^{\pre}_k)$
are $\Sigma$-split.}
 \item{ If $\CA$ is an asymmetric operad then $n\mapsto\CA(n)\otimes\Sigma_n$
is a $\Sigma$-split operad. In particular, the operad for associative algebras
is $\Sigma$-split over $\mathbb{Z}$.}
\end{itemize}

\subsubsection{}
\begin{thm}{cmc-alg}
Let $X$ be a site having enough points. Then any $\Sigma$-split operad 
$\CO\in\Op(C(X^{\pre}_k))$ is homotopically amenable.
\end{thm}

\begin{proof}
It is convenient to sheafify all the picture. If $\CO$ is a presheaf
of $\Sigma$-split dg operads, $\CO^a$ is a sheaf of $\Sigma$-split 
dg operads and $A^a$ is an $\CO^a$-algebra. Sheafification also commutes 
with the free algebra functor and with the coproducts.

We have to check that any generating acyclic cofibration $j:K\to L$
gives rise to a weak equivalence $\alpha$ in the diagram~(\ref{defn:hamo}).
To check this it is sufficient to check that
any fiber functor $\phi: X^{\she}\to\Ens$ transforms the sheafification 
$\alpha^a$ into a quasi-isomorphism.

A fiber functor transforms the diagram~(\ref{defn:hamo}) into a cocartesian
diagram over a ring $k$ with an acyclic complex $K$ and with $L$ being 
the cone of $\id_K$. Note that $K=\phi(K_{\epsilon,n})$ is a 
non-positively graded acyclic complex of $k$-modules. Moreover,
$K$ admits an explicit presentation, see~\cite{sga4}, IV.6.8.3, 
as a filtered direct limit of non-positively graded acyclic complexes 
of free $k$-modules.

Since the free algebra functor commutes with filtered colimits, it is enough
for us to check that the map $\alpha$ in the diagram~(\ref{defn:hamo})
is a weak equivalence provided $K$ is a non-negatively graded complex
of free $k$-modules, $L=\Cone(\id_K)$ and $\CO$ is $\Sigma$-split operad
in $C(k)$. But this latter claim follows from the homotopical amenability 
of $\Sigma$-split operads over $k$.
\end{proof}


In what follows the following definition will be used.

\subsubsection{}
\begin{defn}{sc-sac}A map $f:A\to B$ is called {\em a standard cofibration}
(resp., {\em a standard acyclic cofibration}) if it can be presented as
a direct limit $B=\underset{i\in\mathbb{N}}{\dirlim} A_i$ with $A_0=A$ 
and $A_{i+1}$
being a coproduct of generating cofibrations (resp., of generating
acyclic cofibrations) over $A_i$.  
\end{defn}

We have the following
\subsubsection{}
\begin{prop}{isretract}
Any cofibration is a retract of a standard cofibration.
Any acyclic cofibration is a retract of a standard acyclic cofibration. 
\end{prop}

This, in fact, is true for any cofibrantly generated closed model category.

{\bf In the rest of this section we suppose that the operad $\CO$ is
homotopically amenable.
}

\subsection{Modules}

In this subsection we sketch a presheaf version of~\cite{haha}, Section~5.

\subsubsection{Enveloping algebra}
Enveloping algebra $U(A)$ of an operad algebra $A\in\Alg(\CO,X^{\pre})$
is defined in a usual way. This is an associative algebra in the 
category of complexes $C(X^{\pre}_k)$, such that $U(A)$-modules
are just the modules over the operad algebra $A$.

In particular, the category of modules $\Mod(\CO,A)$ admits
a model structure as in~\ref{O-mod}. Moreover, a presheaf $\cphom_A(M,N)$
is defined for a pair $(M,N)$ of $A$-modules.

\subsubsection{} Weakly equivalent operad algebras have sometimes
non-equivalent derived categories of modules
even for $X^{\pre}=\Ens$, see~\cite{haha}, Section~5.

To get a ``correct'' derived category of modules, one has to work
with cofibrant algebras (or with cofibrant operads and flat algebras,
see~\cite{haha}, 6.8).

\subsubsection{}
\begin{lem}{ac.of.c}
Suppose $\CO$ is $\Sigma$-split.
Let $f:A\to B$ be an acyclic cofibration of cofibrant $\CO$-algebras.
Then the induced map $U(f):\ U(A)\to U(B)$ is a weak equivalence. 
\end{lem}
\begin{proof}The proof is basically the same as that of 
Corollary~5.3.2, \cite{haha}. 
Everything reduces to the case $A$ is generated by a finite number of sections
$x_i$, $i=1,\ldots,n$, over $U_i\in X$ and $B$ is the colimit of a diagram
$$ A\leftarrow F(M)\to F(\id_M),$$
where $M$ is a coproduct of shifts of representable presheaves and has 
acyclic sheafification. The enveloping algebra $U(A)$ admits a filtration
numbered by multi-indices $d:\{1,\ldots, n\}\to\mathbb{N}$, with associated
graded pieces given by the formula
\begin{equation*}
\label{fforA}
\gr_d(U(A))=\CO(|d|+1)\otimes_{\Sigma_d}\bigotimes U_i^{\otimes d_i},
\end{equation*} 
where, as usual, $|d|=\sum d_i$, $\Sigma_d=\prod \Sigma_{d_i}\subseteq 
\Sigma_{|d|}$.

The enveloping algebra $U(B)$ admits a filtration indexed by pairs 
$(d,k)$ with $d$ as above and $k\in\mathbb{N}$. The associated graded piece
corresponding to $(d,k)$ takes form
\begin{equation*}
\label{fforB}
\gr_{d,k}(U(B))=\CO(|d|+k+1)\otimes_{\Sigma_k\times\Sigma_d}
M[1]^{\otimes k}\otimes\bigotimes U_i^{\otimes d_i}.
\end{equation*}
Since $\CO$ is $\Sigma$-split and $M^a$ is flat and acyclic, 
the sheafifications of $ \gr_{d,k}(U(B))$ are acyclic for $k>0$. 
This proves the lemma.
\end{proof}

\subsubsection{}
\begin{cor}{we.of.c}Let $\CO$ be $\Sigma$-split. Any weak equivalence
$f:A\to B$ of cofibrant $\CO$-algebras induces a weak equivalence $U(f)$
of the enveloping algebras.
\end{cor}
\begin{proof}We already know the claim in the case $f$ is an acyclic 
cofibration. Therefore, it suffices to prove it for acyclic fibrations.

Let $f:A\to B$ be an acyclic fibration of cofibrant algebras. There exist
a map $g:B\to A$ splitting $f$: $fg=\id_B$. This implies that for any weak
equivalence $f:A\to B$ the induced map $H(U(A)^a)\to H(U(B)^a)$ 
of the homologies of the sheafifications splits. Applying this to $g:B\to A$
we deduce that it is in fact invertible.  
\end{proof}

\subsection{Differentials and derivations}
\label{diff-der}

In this subsection we present a presheaf version of parts 
of ~\cite{haha}, 7.2, 7.3.

\subsubsection{}
\begin{defn}{der} Let $\CO\in\Op(C(X^{\pre}_k))$ be a presheaf of operads
on $X$, $\alpha:B\to A$ be a map in $\Alg(\CO,X^{\pre})$ and let $M$
be an $A$-module. The presheaf $\cpder_B^{\CO}(A,M)$ of $\CO$-derivations
over $B$ from $A$ to $M$ is defined as the subpresheaf of
$\cphom_k(A,M)$ consisting of local sections which are $\CO$-derivations
from $A$ to $M$ vanishing at $B$.
\end{defn}

By definition $\cpder_B^{\CO}(A,M)$ is a subcomplex of $\cphom_k(A,M)$.
The functor \\
$M\mapsto \cpder_B^{\CO}(A,M)$ is representable in the following 
sense.

\subsubsection{}
\begin{lem}{diff}There exists a (unique up to a unique isomorphism)
$A$-module $\Omega_{A/B}$ together with a global derivation
$\partial:A\to\Omega_{A/B}$ inducing a natural isomorphism in $C(X^{\pre})$
$$ \cphom_A(\Omega_{A/B},M)\isom\cpder_B^{\CO}(A,M).$$
\end{lem} 

The proof of \Lem{diff} is standard, see Proposition 7.2.2 of~\cite{haha}.

The following lemma is the key to the calculation of $\Omega_{A/B}$.

\subsubsection{}
\begin{lem}{as732}
Let $\alpha:B\to A$ be a map of $\CO$-algebras, $M\in C(X^{\pre}_k)$.
Let $f:M\to A$ be a map in $C(X^{\pre}_k)$ and let $A'=A\langle M,f\rangle$
be defined by the cocartesian diagram
$$
\begin{diagram}
 F(\CO,M) & \rTo & F(\CO,\Cone(\id_M)) \\
      \dTo    &               &  \dTo    \\
       A      & \rTo  &  A' 
\end{diagram}
.$$
Put $U=U(\CO,A)$ and $U'=U(\CO,A')$. The map $\partial\circ f:M\to\Omega_{A/B}$
defines $f':U'\otimes M\to U'\otimes_U\Omega_{A/B}$. 

Then the module of differentials $\Omega_{A'/B}$ is naturally isomorphic to 
the cone of $f'$. 
\end{lem}

For the proof see Lemma 7.3.2 of~\cite{haha}. \qed

\subsubsection{}
\begin{prop}{omega(cof)}Let $\CO$ be homotopically amenable.
Let $\alpha:B\to A$ be a cofibration in $\Alg(\CO,X^{\pre})$. Then 
$\Omega_{A/B}$ is cofibrant in $\Mod(\CO,X^{\pre})$.
\end{prop}
\begin{proof}
One can easily reduce the claim to the case when $A$ is generated over $B$ by 
a section $a$ over $U\in X$ subject to a condition $da=b\in B(U)$.
In this case $\cpder_B^{\CO}(A,M)=\cphom_k(U,M)$ so that
$\Omega_{A/B}$ is isomorphic to  $U(A)\otimes U$.
\end{proof}

\subsubsection{}
\begin{cor}{der-fib}
Let $\CO$ be homotopically amenable.
Let $\alpha:B\to A$ is a cofibration in $\Alg(\CO,X^{\pre})$
and let $M$ be a fibrant $A$-module. Then $\cpder_B^{\CO}(A,M)\in 
C(X^{\pre}_k)$ is fibrant.
\end{cor}

\subsubsection{}
Let  $C\overset{\alpha}{\to}B\overset{f}{\to}A$ be a pair of morphisms
in $\Alg(\CO,X^{\pre}_k)$. For each $A$-module $M$ every derivation
$\partial:A\to M$ over $C$ defines a derivation $\partial\circ f:B\to M$
over $C$. This defines a canonical map
\begin{equation*}
\label{omegaf}
\Omega^f: U(B)\otimes_{U(A)}\Omega_{B/C}\to\Omega_{A/C}.
\end{equation*}

The following proposition shows that the module of differentials $\Omega_{A/B}$
has a correct homotopy meaning for cofibrant morphisms $B\to A$. It generalizes
Proposition 7.3.6 of~\cite{haha} where the case $X^{\pre}=\Ens$ is considered. 

\subsubsection{}
\begin{prop}{as736}
Let $C\overset{\alpha}{\to}B\overset{f}{\to}A$ be a pair of morphisms
in $\Alg(\CO,X^{\pre}_k)$. If $f$ is a weak equivalence and 
$\alpha,f\circ\alpha$ are cofibrations than the map $\Omega^f$ is a weak
equivalence.
\end{prop}
\begin{pf}
First of all one proves the claim in the case $f$ is an acyclic cofibration.
This easily follows from~\Lem{as732}. Then one proves the assertion in the 
case the algebras $A$ and $B$ are fibrant. Here the proof of
Proposition 7.3.6 of~\cite{haha} can be repeated verbatim. Finally,
the general case reduces  to the case $A$ and $B$ are fibrant 
passing to fibrant resolutions of $A$ and $B$.
\end{pf}

\subsubsection{}
Let $\CO$ be homotopically amenable.
Let $A$ be a cofibrant $\CO$-algebra on $X$ and let $M$ be 
a fibrant $A$-module. According to \Cor{der-fib},
$\cpder(A,M)\in C(X^{\pre}_k)$ is fibrant. We want to show that 
this object behaves well under localizations.

Let $U\in X$ and let $j:X/U\to X$ be the obvious embedding. We claim that
the presheaf $\cpder(j^*(A),j^*(M))=j^*\cpder(A,M)\in C((X/U)^{\pre}_k)$
has ``the correct homotopy meaning''. As in~\ref{restriction:OK}, the only
problem is that the functor $j^*$ does not always preserve cofibrations.
Therefore, we claim the following.

\begin{lem}{derrestriction-OK}
Let $f:A'\to j^*(A)$ be a cofibrant resolution of $j^*(A)$. Then the 
restriction map
$$ \cpder(j^*(A),j^*(M))\to\cpder(A',j^*(M))$$
is a weak equivalence.
\end{lem}
\begin{proof}
It is enough by ~\Lem{restriction:OK} to check that the map
\begin{equation}
\label{comparingO}
\Omega^f:U(j^*(A))\otimes_{U(A')}\Omega_{A'}\to \Omega_{j^*(A)}
\end{equation}
is a weak equivalence. 

Recall that $A$ is cofibrant. One easily reduces the claim to the case $A$ 
is standard cofibrant. This means that $A$ is presented as a colimit
of $A_n,\ n\in\mathbb{N}$ with $A_0=\CO(0)$ (the initial $\CO$-algebra)
and $A_{n+1}$ defined as the colimit of a diagram
$$ A_n\leftarrow F(\CO,M_n)\to F(\CO,\Cone(\id_{M_n})),$$
$M_n$ being a direct sum of representable presheaves and their shifts.

Then $j^*(A)$ is the colimits of $j^*(A_n)$ with
$j^*(A_{n+1})$ isomorphic to the colimit of the diagram
$$ j^*(A_n)\leftarrow F(\CO,j^*(M_n))\to F(\CO,\Cone(\id_{j^*(M_n)})).$$

If we choose cofibrant resolutions $M'_n\to j^*(M_n)$, one defines recursively
the collection of cofibrant algebras $A'_{n+1}$ as colimits of the diagram
$$ A'_n\leftarrow F(\CO,M'_n)\to F(\CO,\Cone(\id_{M'_n})).$$

Then by induction on $n$ one checks using \Lem{as732} that the map

\begin{equation*}
\label{comparingOn}
\Omega^{f_n}:U(j^*(A_n))\otimes_{U(A'_n)}\Omega_{A'_n}\to \Omega_{j^*(A_n)}
\end{equation*}
is a weak equivalence for each $n$. Passing to a limit, we get a weak 
equivalence~(\ref{comparingO}) for a special choice of $A'=\dirlim A'_n$.

\end{proof}

\subsection{Simplicial structure}
\label{alg-sim}
Similarly to the case $X=*$ described in~\cite{haha}, sect.~4.8, we can define
a simplicial structure on $\Alg(\CO,X^{{\pre}}_k)$ provided 
$k\supseteq\mathbb{Q}$.

\subsubsection{}
\label{simpl-path}
Let $S$ be a finite simplicial set and let $A\in\Alg(\CO,X^{{\pre}}_k)$.
We define the presheaf $A^S$ by the formula
\begin{equation*}
\label{path}
A^S(U)=\Omega(S)\otimes A(U)
\end{equation*}
where $\Omega(S)$ denotes the commutative dg $k$-algebra of polynomial
differential forms on $S$.

\subsubsection{}
For $A,B\in\Alg(\CO,X^{{\pre}}_k)$ define $\shom(A,B)\in\simpl$
by the formula
\begin{equation*}
\label{shom}
\Hom_{\simpl}(S,\shom(A,B))=\Hom(A,B^S) 
\end{equation*}
where $S\in\simpl$ is finite.

The following theorem says that the simplicial structure defined
satisfies Quillen's axiom (SM7).
\subsubsection{}
\begin{thm}{sm7:alg}
Let $\alpha:A\to B$ be a cofibration and $\beta: C\to D$ is a fibration
in $\Alg(\CO,X^{{\pre}}_k)$. Then the natural map
\begin{equation*}
\label{eq:sm7:alg}
\shom(B,C)\to\shom(A,C)\times_{\shom(A,D)}\shom(B,D)
\end{equation*}
is a Kan fibration. It is a weak equivalence if $\alpha$ or $\beta$ is a 
weak equivalence.
\end{thm}

\Thm{sm7:alg} results from the following

\subsubsection{}
\begin{lem}{lemsm7}Let $\alpha:K\to L$ be an injective map of finite 
simplicial sets and let $\beta:C\to D$ be a fibration in 
$\Alg(\CO,X^{{\pre}}_k)$.
Then the natural map
\begin{equation}
\label{eq:lemsm7}
C^L\to C^K\times_{D^K}D^L
\end{equation}
is a fibration. It is weak equivalence if $\alpha$ or $\beta$ is a weak 
equivalence.
\end{lem}
\begin{proof}
The map~(\ref{eq:lemsm7}) can be rewritten as 
$$ \Omega(L)\otimes C\to \Omega(K)\otimes C
\times_{\Omega(K)\otimes D} \Omega(L)\otimes D.$$
it is pointwise surjective since $\beta$ is pointwise surjective.
This also implies that~(\ref{eq:lemsm7}) is an acyclic fibration provided 
$\alpha$ or $\beta$ is a weak eqivalence. 

Since ~(\ref{eq:lemsm7}) is pointwise surjective, it is enough
to check that its kernel is fibrant. The kernel easily identifies
with the tensor product $K_{\alpha}\otimes K_{\beta}$ where 
$K_{\alpha}=\Ker\{\alpha:\Omega(L)\to\Omega(K)\}$ and $K_{\beta}=\Ker(\beta)$.

Presheaf $K_{\beta}$ is fibrant. Complex of $k$-modules $K_{\alpha}$
has form $k\otimes_{\mathbb{Q}}K^{\mathbb{Q}}_{\alpha}$ where
the complex $ K^{\mathbb{Q}}_{\alpha}$ 
has finite dimensional cohomology. One can therefore
write $K_{\alpha}=H\oplus C$ where $H$ is the finite dimensional cohomology
of $K_{\alpha}$ and $C$ is contractible. Then $H\otimes K_{\beta}$ is fibrant
as a finite direct sum of (shifts of) $K_{\alpha}$; $C\otimes K_{\beta}$
is pointwise contractible, and therefore fibrant.
\end{proof}

\subsubsection{Functor $\Tot$}
Let $A^{\bullet}$ be a cosimplicial object in $\Alg(\CO,X^{\pre})$.
The algebra $\Tot(A^{\bullet})$ is defined by the standard formula
\begin{equation*}
\label{tot:alg}
\Tot(A^{\bullet})=\invlim_{p\to q}(A^q)^{\Delta^p}.
\end{equation*}

\subsection{Descent}
\subsubsection{}
 Let $V_{\bullet}$ be a hypercover of
$X$ and let $\CO$ be an operad in $X^{\she}_k$. Put
$\CO_n=\CO|_{V_n}$. Let $\Alg^{\f}$ (resp., $\Alg^{\f}_n$) denote the 
category of fibrant $\CO$-algebras on $X$ (resp., fibrant
$\CO_n$-algebras on $X/V_n$). 

The assignment 
$n\longmapsto\Alg^{\f}_n$ defines a
cosimplicial object in $\Cat$ (more precisely, a category cofibrant 
over $\Delta$). We denote by $\Alg^{\f}(V_{\bullet})$ the following category.
The objects of $\Alg^{\f}(V_{\bullet})$ are collections 
$\{A_n\in\Alg^{\f}_n\}$
together with weak equivalences $A_n\to\phi^*(A_m)$ corresponding to each
$\phi:m\to n$ in $\Delta$, satisfying the standard cocycle condition.

The functor $\epsilon^*:\Alg^{\f}\to\Alg^{\f}(V_{\bullet})$ assigns to
each algebra on $X$ the collection of its restrictions. This functor
preserves weak equivalences. We define the functor
\begin{equation*}
\label{alpha_*}
\alpha_*:\Alg^{\f}(V_{\bullet})\to\Alg^{\f}
\end{equation*}
by the formula
\begin{equation*}
\label{alpha_*=}
\alpha_*(A_n,\phi)=\Tot\{n\mapsto \iota^n_*(A_n)\}.
\end{equation*}  

Here the functor $\Tot$ is defined as in~\ref{alg-sim} and $\iota^n_*$ is the
direct image of the localization functor $\iota^n: X/V_n\to X$, 
\cite{sga4}, III.5.

\subsubsection{}
\begin{prop}{prop:desc}The functors $\alpha_*$ and $\alpha^*$ induce an 
adjoint pair of equivalences on
the corresponding Dwyer-Kan localizations of $\Alg^{\f}(V_{\bullet})$ and
 $\Alg^{\f}$ with respect to quasi-isomorphisms. 
\end{prop}
\begin{proof}
The functors involved do not depend on the operad $\CO$. Therefore, it is 
enough to prove the claim for the trivial operad so that $\Alg(\CO)$
is just the category of complexes. This is done in~\cite{sga4}, Expos\'e Vbis,
for $C^+(X^{\she}_k)$ and in~\cite{hs}, Sect.~21, for unbounded complexes.

\end{proof}

\section{Deformations}
\label{sect:result}

Let $k$ be a field. In this section we define a functor of formal
deformations of a sheaf of operad algebras over $k$. In the case
when $X$ has enough points and $k\supseteq\mathbb{Q}$ we present a
dg Lie algebra which governs (under some extra restrictions on the sheaf of 
algebras $A$) the described deformation problem.

\subsection{Deformation functor}

\subsubsection{}
Let $X$ be a site, $\CO\in\Op(C(X^{\she}_k))$ be a sheaf of dg operads on $X$
and let $A\in\Alg(\CO,X^{\she})$ be a sheaf of $\CO$-algebras.
The functor we define below describes formal deformations of the sheaf
of $\CO$-algebras $A$.

\subsubsection{Bases of deformations}
Fix a commutative algebra homomorphism $\pi:K\to k$ (it will usually
be the identity $\id:k\to k$ of a field of characteristic zero).

Let $\dgar(\pi)$ denote the category of commutative non-positively graded dg
algebras $R$ endowed with morphisms $\alpha:K\to R,\ \beta:R\to k$, such
that $\Ker(\beta)$ is a nilpotent ideal of finite length.
The category $\dgar(\pi)$ is the category of allowable bases for formal
deformations of $A$.

We suppose that the operad $\CO$ is obtained by the base change from a fixed
operad $\CO_K\in\Op(C(X^{\she}_K))$, so that $\CO=\CO_K\otimes k$.
For each $R\in\dgar(\pi)$ we put $\CO_R=\CO_K\otimes R$.

\subsubsection{}
Fix $R\in\dgar(\pi)$. Let $\Alg^{\fl}(\CO_R,X^{\she})$ denote the
category of sheaves of $\CO_R$-algebras {\em flat as sheaves of
$R$-modules}. Let $\CW_{\she}^{\fl}(R,X)$ be the subcategory of
weak equivalences in $\Alg^{\fl}(\CO_K\otimes R,X^{\she})$ (i.e.
quasi-isomorphisms of complexes of sheaves). In what follows we
will usually omit $X$ from the notation.

Recall that Dwyer-Kan construction~\cite{dk,dk2,dk3} assigns to a pair
$(\CC,\CW)$ where $\CW$ is a subcategory of $\CC$, a {\em hammock
localization} $L^H(\CC,\CW)$ which is a simplicial category. We will apply
this construction to pairs of coinciding categories $(\CW,\CW)$. The
Dwyer-Kan construction gives in this case a simplicial groupoid.
We will denote it $\widehat{\CW}$ and call {\em a weak groupoid completion},
see Appendix~\ref{App:scat} for details.

The base change functor
$$ M\mapsto M\otimes_Rk$$
induces a functor between the weak groupoid completions
\begin{equation}
\label{base-change}
\beta^*: \widehat{\CW}_{\she}^{\fl}(R)\to\widehat{\CW}_{\she}(k).
\end{equation}

\subsubsection{}
\begin{defn}{def-functor}
Deformation functor of $\CO$-algebra $A$ is the functor
\begin{equation*}
\label{eq:def-functor}
\cdef_A:\dgar(\pi)\to\sGrp
\end{equation*}
from the category $\dgar(\pi)$ to the category of simplicial groupoids $\sGrp$
defined as the homotopy fiber of~(\ref{base-change}) at
$A\in\widehat{\CW}_{\she}(k).$
\end{defn}
Since~(\ref{base-change}) is a functor between fibrant simplicial categories,
the homotopy fiber is represented by the fiber product
\begin{equation*}
\label{def-functor-presentation}
*\times_{\widehat{\CW}_{\she}(k)}\widehat{\CW}_{\she}(k)^{\Delta^1}
\times_{\widehat{\CW}_{\she}(k)}\widehat{\CW}_{\she}^{\fl}(R),
\end{equation*}
see~\ref{hfib.scat}. This allows one to consider $\cdef_A$ as a functor
with values in $\sGrp$ and not in the corresponding homotopy category.

\subsection{Properties of $\cdef_A$}

We list below some properties of the functor $\cdef_A$ which justify
the above defintion.

\subsubsection{}
The functor $\cdef_A$ does not depend on the specific choice
of the site $X$; only the topos $X^{\she}$ counts.

\subsubsection{}
Suppose that $X$ is a one-point site, so that we are dealing with deformations
of an operad algebra $A$. Suppose also that $k$ is a field of characteristic
zero and $\pi=\id_k:K=k\to k$.

Recall that in~\cite{dha} a deformation functor is defined for this setup.
In the definition the simplicial category of cofibrant algebras
and weak equivalences is used instead of the hammock localization.
We prove in~\Prop{eq-weak-gr} that these two functors are equivalent.

\subsubsection{Descent} 
\label{def.descent}
Let $V_{\bullet}$ be a hypercover of
$X$ and let $A$ be a sheaf of operad algebras on $X$. Put
$A_n=A|_{V_n}$. The assignment $n\longmapsto\cdef_{A_n}$ defines a
cosimplicial object in $\sCat$. One has a natural descent functor
\begin{equation*}
\label{eq:descent} 
\cdef_A\to\holim\{n\mapsto\cdef_{A_n}\},
\end{equation*}
see~\ref{holim.scat}.
We claim the functor described is an equivalence. 
In fact, since
$\cdef_A$ is defined as a homotopy fiber, the right-hand side
of~(\ref{eq:descent}) is the homotopy fiber of the map
\begin{equation*}
\label{}
\holim\{n\mapsto\widehat{\CW}_{\she}^{\fl}(R,V_n)\}\to
\holim\{n\mapsto\widehat{\CW}_{\she}^{\fl}(k,V_n)\}.
\end{equation*}

It is enough therefore to check that the functor 
\begin{equation*}
\label{}
\widehat{\CW}_{\she}^{\fl}(R,X)\to
\holim\{n\mapsto\widehat{\CW}_{\she}^{\fl}(R,V_n)\}
\end{equation*}
is an equivalence. Since all simplicial categories involved are simplicial
groupoids, it is enough by~\ref{eq-lcc-ldd} to check that this functor 
induces an equivalence of the nerves. 
The functor $\holim$ commutes with the nerve functor,
see~\cite{uc}, Prop. A.5.2, or~\ref{holim-vs-N}. 
Moreover, the nerve of a Bousfield-Kan localization is equivalent 
to the nerve of an original category. Therefore, the claim 
follows from~\Prop{prop:desc} and \ref{eq-weak-gr}.

\subsubsection{Connected components}

Here we assume that the topos $X^{\she}$ admits enough points.
We assume as well that $k$ is a field of characteristic zero and
$\pi=\id_k$.

Suppose that $R$, $A$ and $\CO$ are concentrated at degree zero, so that
we have a classical deformation problem. We claim that $\Def_A(R)$
is equivalent to the groupoid $\Def^{\cl}_A(R)$ of flat $R$-deformations 
of $A$.

Let $B$ be a $R$-flat  sheaf of (dg) $R\otimes\CO$-algebras such that
the reduction $\beta^*(B)$ is quasi-isomorphic to $A$. We claim
that $B$ is concentrated in degree zero and $H^0(B)$ is a flat
deformation of $A$. In fact, since $X^{\she}$ has enough points, the
claim can be verified fiberwise. The case $X^{\she}=\Ens$ is explained
in~\cite{dha}, 5.1.

The assignment $B\mapsto H^0(B)$ defines, therefore, a simplicial
functor $\Def_A(R)\to \Def^{\cl}_A(R)$. We claim this functor
is an equivalence. In fact, it is enough by~\ref{eq-lcc-ldd} to check the 
functor induces an equivalence of the nerves. 
Let $\CW_{\she}^{\fl}(R)^0$ denote the full subcategory of 
$\CW_{\she}^{\fl}(R)$ consisting of algebras whose cohomology is flat and
concentrated in degree zero , and let $\Alg^{\iso}(R)$ denote the groupoid of
flat $R\otimes \CO$-algebras concentrated in degree zero.

It is enough to check that the  functor
$$ H^0:\CW_{\she}^{\fl}(R)^0\to\Alg^{\iso}(R)$$
induces an equivalence of the nerves.

To prove this, consider a third category $\CW^{\c}(R)$ consisting
of cofibrant {\em presheaves} $P$ of $R\otimes\CO$-algebras. 
Sheafification defines a functor
$$ a: \CW^{\c}(R)\to\CW_{\she}^{\fl}(R).$$
We denote by $\CW^{\c}(R)^{0,\fl}$ the full subcategory of $\CW^{\c}(R)$
consisting of complexes of presheaves whose sheafification belongs to
$\CW_{\she}^{\fl}(R)^0$. The restriction defines the finctors
$$ a: \CW^{\c}(R)^{0,\fl}\to\CW_{\she}^{\fl}(R).$$
and
$$ H^0\circ a: \CW^{\c}(R)^{0,\fl}\to\Alg^{\iso}(R).$$
Both functors $a$ and $H^0\circ a$ induce an equivalence of nerves 
by Quillen's Theorem A~\cite{hakt1} and \Thm{res-cont}. 

This proves the assertion.

\subsection{Reformulation in terms of presheaves}

From now on we assume that $k$ is a field of characteristic zero and
$\pi=\id_k$. We assume as well that the topos $X^{\she}$ admits enough 
points.

\subsubsection{Notation}
The category $\Alg(\CO\otimes R, X^{\pre})$ admits a
model structure defined in~\ref{cmc-alg:hamo}. We denote the subcategory
of weak equivalences by $\CW(R)$. The notation $\CW^{\c}(R)$, $\CW^{\f}(R)$,
$\CW^{\cf}(R)$ and $\CW^{\cf}_*(R)$ has the same meaning as
in~\ref{eq-weak-gr}.

\Prop{wgs2} below claims that the weak groupoid
$\widehat{\CW}_{\she}^{\fl}(R)$ admits an equivalent description in terms of
the weak groupoid of the category of presheaves. This description
is functorial in $R$. This allows one to conveniently describe the homotopy
fibre of~(\ref{base-change}).

Note the following technical lemma.

\subsubsection{}
\begin{lem}{cofa}Let $A$ be a  cofibrant 
$\CO\otimes R$-algebra and let $A^{\#}$ be the corresponding 
presheaf of $R$-modules. There is an increasing filtration 
$\{F^iA,i\in I\}$ of $A^{\#}$ indexed by a well-ordered set $I$ such that 
the associated graded factors $\gr^i(A)=F^i(A)/\sum_{j<i}F^jA$ are isomorphic
to
\begin{equation*}
  \label{eq:factors}
  R\otimes X^i
\end{equation*}
with $X^i\in C(X^{\pre}_k)$.
\end{lem}
\begin{proof}
We can assume that $A$ is a direct limit over a well ordered set
of generating acyclic cofibrations. In fact, a general cofibrant algebra
is a retract of an algebra of this type, and the property claimed in the
lemma is closed under retractions.

Let, therefore, the algebra $A$ be freely generated
by a collection of sections $x_j$ over $U_j\in X$ of degree $d_j$
numbered by a well-ordered set $J$. We define $I$ to be the collection
of functions $m:J\to\mathbb{N}$ nonzero at a finite number of elements of $J$.
The set $I$ is endowed with the lexicographic order:

\begin{equation*}
m < m'\text{ iff }\exists k\in J: m(k)<m'(k)\text{ and } m(j)=m'(j)
\text{ for }j>k.  
\end{equation*}

The increasing  filtration is numbered by elements of $I$ 
and the associated graded factors have the form
$$R\otimes\CO(n)\otimes_{\Sigma_m}k[-d_1]\cdot U_1^{\otimes m_1}
\otimes\ldots\otimes k[-d_k]\cdot U_k^{\otimes m_k}$$
where $n=\sum m_i,\quad \Sigma_m=\prod\Sigma_{m_i}\subseteq\Sigma_n$,
$\Sigma_n$ being the symmetric group. 
\end{proof}

\subsubsection{}
\begin{lem}{beta^*}
1. The functor $\beta^*:\Alg(\CO_R,X^{\pre})\to\Alg(\CO,X^{\pre})$ preserves
cofibrations.

2. Let $A\in\Alg^{\c}(\CO_R,X^{\pre})$. Then $A$  is fibrant if and only if 
$\beta^*(A)$is fibrant.
\end{lem}

\begin{proof}
1. The first claim is obvious for generating cofibrations. 
The general case follows from the fact that $\beta^*$ commutes with 
colimits and retracts.

2. According to~\Lem{cofa}, $A^{\#}$ admits an increasing filtration 
$\{F^iA,i\in I\}$ with associated graded factors having a form
\begin{equation*}
  \gr^iA=R\otimes X^i.
\end{equation*} 
This implies that $A(U)$ is a cofibrant $R$-module for each $U\in X$.

$A$ is fibrant iff for any hypercover
$\epsilon:V_{\bullet}\to U$ in $X$ the natural map
\begin{equation}
\label{tocech}
A(U)\to\Cech(V_{\bullet},A)
\end{equation}
is a weak equivalence. 
In \Lem{cech-is-cof} below we prove that $\Cech(V_{\bullet},A)$ is a 
cofibrant $R$-module.
The Cech complex of $\beta^*(A)$ is just the reduction of the Cech complex
of $A$ modulo the maximal ideal of $R$. Therefore, if~(\ref{tocech}) is a weak equivalence, the reduction 
modulo the maximal ideal of $R$ is also a weak equivalence. In the other
direction, there exists a finite filtration of $R$ by
$k$-subcomplexes such that the associated graded factors are isomorphic 
to $k$ up to shift. Therefore, the cone of~(\ref{tocech}) admits a 
finite filtration with acyclic associated graded factors.
\end{proof}

\subsubsection{}
\begin{lem}{cech-is-cof}
Let $A\in\Mod(R,X^{\pre}_k)$ admit an increasing filtration 
$\{F^iA, i\in I\}$ with associated graded factors of form 
\begin{eqnarray*}
  \gr^i(A)=R\otimes X^i,
\end{eqnarray*}
for some $X^i\in C(X^{\pre}_k)$. Then for each hypercover 
$\epsilon:V_{\bullet}\to U$ the Cech complex $\Cech(V_{\bullet},A)$
is a cofibrant $R$-module.
\end{lem}
\begin{proof}
Cech complex of $A$ is the total complex corresponding to the
bicomplex 
$$A_0\overset{\delta^0}{\to} A_1\overset{\delta^1}{\to}\ldots$$
where $A_n$ is the collection of sections of $A(V_n)$ vanishing on the 
degenerate part of $V_n$ (see~(\ref{cech}) for the precise definition).

Each $A_n$ is filtered and $\delta^n$ preserves the filtration.
Define $B_n=\ker(\delta^n)$. $B_n$ are filtered $R$-modules with associated
graded factor having the form
$$ R\otimes Y^i_n,\quad Y^i_n=\ker(\delta^n:X^i_n\to X^i_{n+1}).$$
The restriction of the differential $\delta_n$ on $B_n$ vanishing,
the corresponding total subcomplex is isomorphic to
$\prod_nB_n[-n]$. Let us show it is a cofibrant $R$-module.
In fact, each $B_n$ admits an increasing filtration with $R$-free
associated graded factors. Since $R$ is artinian, a product of free $R$-modules
is free. Therefore, $\prod_nB_n[-n]$ admits as well an increasing filtration
with $R$-free associated graded factors. Therefore, it is $R$-cofibrant.
   
The quotient is the total complex of the bicomplex
$$A_0/B_0\to A_1/B_1\to\ldots$$
which has also a vanishing horisontal differential. Each quotient 
$A_n/B_n$ admits a filtration with the associated graded factors
having form $R\otimes X^i_n/Y^i_n$ which is also $R$-cofibrant.
Lemma is proven.

\end{proof}

The following proposition claims that the simplicial categories
$\CW^{\cf}_*(R)$ and of $\widehat{\CW}^{\fl}_{\she}(R)$ are 
canonically equivalent.

\subsubsection{}
\begin{prop}{wgs2}
There is a canonical in $R$  collection of equivalences
of weak groupoids
\begin{equation*}
\label{eq:topre}
\CW^{\cf}_*(R)\overset{c}{\to}\widehat{\CW}^{\cf}_*(R)
\overset{b}{\leftarrow}\widehat{\CW}^{\cf}(R)
\overset{a}{\to}\widehat{\CW}_{\she}^{\fl}(R).
\end{equation*}
\end{prop}
\begin{proof}
The map $c$ is a canonical map from a simplicial category to its hammock 
localization. It is an equivalence since $\CW^{\cf}_*(R)$ is a weak groupoid,
see~\ref{eq-vs-nerves}. The map $b$ is equivalence by part (ii)
of~\Prop{eq-weak-gr}. The map $a$ is induced by the sheafification.

First of all, by~\Lem{cofa} the sheafification of $A\in {\CW}^{c}(R)$
is $R$-flat. Therefore, the functor $a$ is defined. To prove $a$ is an 
equivalence, it suffices, by~\Cor{eq-lcc-ldd}, to check that the functor
$$a:\CW^{\cf}(R)\to\CW^{\fl}_{\she}(R)$$
induces an equivalence of the nerves.

\Prop{eq-weak-gr} (2) asserts that the map
\begin{equation*}
\CW^{\cf}(R)\to \CW^{\c}(R)
\end{equation*}
induces an equivalence of the nerves.

The sheafification functor $a: \CW^{\c}(R)\to\CW_{\she}(R)$
has its image in $\CW^{\fl}_{\she}(R)$. The resulting
functor
$$a: \CW^{\c}(R)\to\CW^{\fl}_{\she}(R)$$
induces an equivalence of the nerves by~\Thm{res-cont} and Quillen's
Theorem A, see~\cite{hakt1}.

\end{proof}

\subsection{Fibration lemma}

\subsubsection{}
\label{more.restrictions}
In this subsection we assume the following properties.
\begin{itemize}
  \item[$\bullet$]{For each $n\in\mathbb{N}$ the complex of sheaves $\CO(n)$ 
is non-positively graded,  $\CO(n)\in C^{\leq 0}(X^{\she})$.
}
  \item[$\bullet$]{The site $X$ admits a final object.}

\end{itemize}

Let $\CW^{\cf}_*(R)_{\leq 0}$ denote the full simplicial
subcategory of $\CW^{\cf}_*(R)$ consisting of algebras $A$ satisfying
the extra condition $\CH^i(A)=0$  for $i>0$.

\subsubsection{}
\begin{lem}{is-fibration}Suppose that the condition
of~\ref{more.restrictions} on $\CO$ and $X$ are fulfilled.
Then the functor
\begin{equation}
\label{eq:is-fibration}
\beta^*:\CW^{\cf}_*(R)_{\leq 0}\to\CW^{\cf}_*(k)_{\leq 0}
\end{equation}
 is a fibration in $\sCat$.
\end{lem}

We recall the model category structure on $\sCat$ in~\ref{ss:cmc-scat}.
The proof of~\Lem{is-fibration} is similar to that of
Lemma 4.2.1 of~\cite{dha}. It is presented
in~\ref{step1} --~\ref{eopf:isfibration} below.
The properties of derivations studied in~\ref{diff-der} play an
important role in the proof.

\subsubsection{}
\label{step1}
Let $B$ be a cofibrant  presheaf of $\CO$-algebras.
We claim that $B$ and $\alpha^*\beta^*(B)$ are isomorphic as
presheaves of graded (without differential) algebras. In fact, any cofibrant
presheaf is a retract of a standard cofibrant (see~\ref{sc-sac}) presheaf of 
algebras, so we can assume that $B$ is standard. 
This means, in particular, that $B$ is freely generated, as a 
graded presheaf of algebras, by a collection of generators $b_i$ which 
are sections of $B$ over some objects  $U_i\in X$.
In this case the claim is obvious.

Therefore, any cofibrant presheaf $B$ of $\CO$-algebras is isomorphic to
a presheaf of form $(\alpha^*(A), d+z)$ where $(A,d)=\beta^*(B)$ and
$z\in\fm\otimes\der(A,A)$ satisfies the Maurer-Cartan
equation. Here $\fm=\Ker(\beta:R\to k)$ is the maximal ideal of $R$
and $\der$ denotes the global sections of the presheaf $\cpder$.

\subsubsection{}
\label{step2}
A morphism of simplicial categories is a fibration if it satisfies the
conditions (1), (2)  of~\Defn{scat-fib}.
Let us check the condition (1).

Let $f:A\to B$ be a weak equivalence
of fibrant cofibrant algebras in $\Alg(\CO,X^{\pre})$.
Let $\CT_A=\cpder(A,A)$ and similarly for $\CT_B$. The global sections
of the derivation algebras are denoted $T_A$ and $T_B$ correspondingly.
Let one of two elements $a\in\MC(\fm\otimes T_A),\ b\in\MC(\fm\otimes T_B)$
be given. We have to check that there exists a choice of the second element
and a map
\begin{equation}
\label{mapg}
g:(\alpha^*(A),d+a)\to (\alpha^*(B),d+b)
\end{equation}
in $\Alg(R\otimes\CO, X^{\pre})$ lifting $f$.

Note that under the restrictions of~\ref{more.restrictions} any algebra
$(\alpha^*(A),d+a)$ belongs to $\CW^{\cf}(R)$.

We can consider separately the cases when $f$ is an acyclic fibration or an
acyclic cofibration. In both cases we will be looking for the map~(\ref{mapg})
in the form
\begin{equation*}
\label{g=}
g=\gamma_B^{-1}\circ \alpha^*(f)\circ\gamma_A
\end{equation*}
where $\gamma_A\in\exp(\fm\otimes T_A)^0$ and similarly for $\gamma_B$.
A map~(\ref{mapg}) should commute with the differentials $d+a$ and $d+b$.
This amounts to the condition
\begin{equation*}
\label{cond.g}
f_*(\gamma_A(a))=f^*(\gamma_B(b)),
\end{equation*}
where the natural maps
\begin{equation}
\label{nat.maps.for.T}
T_A\overset{f_*}{\lra}\der_f(A,B)\overset{f^*}{\lla} T_B
\end{equation}
are defined as the global sections of the standard maps
\begin{equation*}
\label{nat.maps.for.CT}
\CT_A=\cphom(\Omega_A,A)\overset{f_*}{\lra}\cphom(\Omega_A,B)\overset{f^*}
{\lla}\cphom(\Omega_B,B)= \CT_B.
\end{equation*}
The maps $f_*$ and $f^*$ in~(\ref{nat.maps.for.T}) are weak equivalences as
global sections of weak equivalences between the fibrant presheaves
$\CT_A,\,\CT_B,\,\cpder_f(A,B)$.

Recall that we assume that $f$ is either acyclic cofibration or
acyclic fibration.

\subsubsection{}
\begin{lem}{diagr_exists}
Let $f:A\to B$ be a weak equivalence of fibrant cofibrant algebras. 
Suppose that
$X$ admits a final object.
Suppose that either
\begin{itemize}
  \item[(af)]{ $f$ is an acyclic fibration}
\end{itemize}
or
\begin{itemize}
  \item[(ac)]{ $f$ is a standard acyclic cofibration.}
\end{itemize}
Then there exists a commutative
square

\begin{equation*}
\label{diagr_dgl}
\begin{diagram}
&  &  T_f &  &   \\
& \ldTo^{g} & & \rdTo^{h}& \\
T_A & \rTo^{f_*} & \der_f(A,B) & \lTo^{f^*} & T_B
\end{diagram}
\end{equation*}
where $T_f$ is a dg Lie algebra and $g,\! h$ are Lie algebra
quasi-isomorphisms.
\end{lem}
\begin{proof}
Note first of all that the maps $f^*,\, f_*$ are weak equivalences.
This follows, as for the ``absolute'' case of~\cite{haha}, 8.1, from
the presentation
\begin{equation*}
\label{t-thru-omegas}
T_A=\chom(\Omega_A,A);\quad T_B=\chom(\Omega_B,B);\quad \der_f(A,B)=
\chom(\Omega_A,B).
\end{equation*}

We construct the dg Lie algebra $T_f$ as follows.

{\em Case 1.} $f$ is acyclic fibration. Let $I=\Ker(f)$. Define
$T_f$ to be the subalgebra of $T_A$ consisting of (global) derivations
preserving $I$. Since $A$ is cofibrant, any global derivation of $B$
can be lifted to $A$. Therefore, the map $h:T_f\to T_B$ is surjective.
The kernel of $h$ consists of global derivations on $A$ with values in $I$.
The presheaf $\cpder(A,I)=\shom(\Omega_A,I)$ is acyclic fibrant. Therefore,
its global sections are acyclic since $X$ admits a final object.

{\em Case 2.} $f$ is a standard acyclic cofibration.
We define $T_f$ as the collection of
(global) derivations of $B$ preserving $f(A)$. Let us check that the map
$g:T_f\to T_A$ is a surjective quasi-isomorphism.
The algebra $B$ is obtained from $A$ by a sequence of generating acyclic
cofibrations corresponding to  hypercovers of $X$.

Let $B=\dirlim B_i$ with $B_0=A$ and $B_i$ obtained from the cocartesian
diagram
\begin{diagram}
F(K_i) & \rTo^{F(j_i)} & F(L_i) \\
\dTo^{F(\phi_i)} & & \dTo  \\
B_{i-1} & \rTo & B_i,
\end{diagram}
where $j_i:K_i\to L_i$ are acyclic cofibrations of presheaves and $F(\quad)$
is the free algebra functor.

Put $T_i=\{\delta\in\der(B_i,B)|\delta(f(A))\subseteq f(A)\}$.
Then $T_0=T_A$ and $T_f=\invlim T_i$. We claim that the maps
$g_i:T_i\to T_{i-1}$ are acyclic fibrations. Then $g:T_f\to T_A$
is also an acyclic fibration.

To prove surjectivity of $g_i$ fix a derivation $\delta: B_{i-1}\to B$.
 Derivations of $B_i$ extending $\delta$
correspond to commutative diagrams of presheaves
\begin{center}
$$
  \begin{diagram}
    K_i &  \rTo^{j_i}  & L_i \\
   \dTo^{\phi_i} & & \dTo \\
    B_{i-1} & \rTo^{\delta} & B
  \end{diagram}
$$
\end{center}
Since the map $j_i$ is an acyclic cofibration and $B$ is fibrant, there
is no obstruction to extending $\delta$. Surjectivity of $g_i$ is proven.

Now the kernel of $g_i$ identifies with $\chom(L_i/K_i,B)$. The presheaf
$\cphom(L_i/K_i,B)$ being acyclic and fibrant, its global sections
$\chom(N/M,B)$ are acyclic.
\Lem{diagr_exists} is proven.
\end{proof}

Recall now (see~\ref{deligne}) that for a dg Lie algebra $\fg$ the 
formal groupoid
\begin{equation*}
\label{deligne_gr}
\Del_{\fg}:\dgar(k)\to\Grp
\end{equation*}
is defined as the transformation groupoid of the group $\exp(\fm\otimes\fg)^0$
acting on the set of Maurer-Cartan elements of $(\fm\otimes\fg)^1$.

The set of connected components $\pi_0(\Del_{\fg}(R))$ is a weak homotopy
invariant of $\fg$. Thus, the maps $g$ and $h$ induce bijections
$$ \pi_0(\Del_{T_A}(R))\lla\pi_0(\Del_{T_f}(R))\lra\pi_0(\Del_{T_B}(R)).
$$
of the sets of components. This proves the condition (1)
of~\ref{scat-fib}.

\subsubsection{}
\label{eopf:isfibration}
Let us check the condition (2) of~\ref{scat-fib}.

Let $\widetilde{A},\,\widetilde{B}\in\CW^{\cf}(R)$ and let
$A=\beta^*(\widetilde{A}),\ B=\beta^*(\widetilde{B})$. We have to check that
the map
\begin{equation}
\label{reduction}
\shom(\widetilde{A},\widetilde{B})\to\shom(A,B)
\end{equation}
is a Kan fibration. The algebra $B$ can be considered as $\CO_R$-algebra
with $R$ acting on $B$ through $\beta:R\to k$. The canonical
map $\widetilde{B}\to B$ is pointwise surjective
and both $B$ and $\widetilde{B}$ are fibrant. Therefore, it is
a fibration by~\Thm{gen-com} (2). Since $\widetilde{A}$ is cofibrant,
the map
\begin{equation}
\label{reduction2}
\shom(\widetilde{A},\widetilde{B})\to\shom(\widetilde{A},B)
\end{equation}
is a Kan fibration. But the maps~(\ref{reduction}) and~(\ref{reduction2})
coincide. This proves the condition (2) of~\ref{scat-fib}.

Fibration \Lem{is-fibration} is proven.
\qed


\subsection{The main theorem.}

\subsubsection{Assumptions}
\label{main-assumptions}
In this subsection we assume that the following conditions.

\begin{itemize}

\item{For each $n\in\mathbb{N}$ the complex $\CO(n)$ of sheaves
is non-positively graded,  $\CO(n)\in C^{\leq 0}(X^{\she})$. }

\item{ The $\CO$-algebra $A$ satisfies the property $\CH^i(A)=0$ for $i>0$. }

\end{itemize}

\subsubsection{}
\begin{lem}{fc-repr}Under the assumptions of~\ref{main-assumptions},
there exists an algebra $A'$ weakly equivalent to $A$ such that
\begin{itemize}
  \item{$A'$ is fibrant and cofibrant;}
  \item{$A'\in C^{\leq 0}(X^{\pre}_k)$.}
\end{itemize}
\end{lem}
\begin{proof}
Let $A\to B$ be a fibrant resolution of $A$. One has $H^i(B)=0$ for $i>0$.
Let $C=\tau^{\leq 0}(B)$. The canonical map $C\to B$ is therefore pointwise
quasi-isomorphism. Then \Thm{gen-com} (2) asserts that $C$ is fibrant as well.
One can  choose a cofibrant resolution $A'$ of $C$ having generators
only in non-positive degrees. The algebra $A'$ satisfies all necessary
properties.
\end{proof}

\subsubsection{}Let $A$ be as above. Let $A'$ be an algebra whose existence
is guaranteed by \Lem{fc-repr}. We define the {\em local} tangent Lie
algebra of $A$ by the formula
\begin{equation*}
\label{ltla}
\CT_A=\cpder(A',A').
\end{equation*}

This is a fibrant presheaf of dg Lie algebras. Define, finally, {\em global}
tangent Lie algebra of $A$ by the formula

\begin{equation*}
\label{gtla}
T_A=\Right\Gamma\CT_A
\end{equation*}
where the functor $\Right\Gamma$ assigning a dg Lie algebra to a presheaf
of dg Lie algebras, is understood in the following sense (compare 
to~\ref{RGamma}).
Choose a hypercover $V_{\bullet}$ of $X$. The collection
\begin{equation*}
\label{global-t}
n\mapsto \Hom(V_n,\CT_A)
\end{equation*}
is a cosimplicial dg Lie algebra. The corresponding $\Tot$ functor
produces a dg Lie algebra which is denoted $\Right\Gamma\CT_A$.
The result does not depend on the choice of $V_{\bullet}$.

Note that if the site $X$ admits a final object $*$, $T_A$ is equivalent
to  $\CT_A(*)$.

\subsubsection{}

A simplicial presheaf $K_{\bullet}$ is called {\em finite dimensional}
if it coincides with its $n$-th skeleton  for some $n$. 

In \Thm{result} below we require the site $X$  admit a finite
dimensional hypercover. This condition is void if $X$ admits
a final object. In the case $X$ is the site of affine open subschemes
of a scheme  $S$, the requirement is fulfilled if $S$ is quasi-compact 
separated or finite dimensional scheme.



\subsubsection{}
\begin{thm}{result}
Suppose that the conditions of~\ref{main-assumptions} on $\CO$ and $A$ are
fulfilled. Suppose also that the site $X$ admits a finite dimensional
hypercover.

Then the deformation functor $\cdef_A$ is equivalent to
$\cdel_{T_A}$ where $T_A$ is the (global) tangent Lie algebra of $A$.
\end{thm}

\subsubsection{}
\begin{rem}{result.for.final}
Fibration Lemma~\ref{is-fibration} implies \Thm{result} if $X$ admits
a final object. In fact, let $A$ be a fibrant cofibrant $\CO$-algebra
with $A^i=0$ for $i>0$. By~\Lem{is-fibration} the homotopy fiber 
of~(\ref{eq:is-fibration}) is equivalent to its usual fiber.
Fix $(R,\fm)\in\dgar(k)$. The fiber of~(\ref{eq:is-fibration}) at $A$
is the simplicial groupoid defined by the perturbations of the differential
in $R\otimes A$. This is precisely the simplicial Deligne groupoid
$\cdel_T(R)$ defined in~\cite{dha} (see~\ref{deligne}), with $T=\der(A,A)$.
\end{rem}

\begin{proof}[Proof of the theorem]

In a few words, the proof is the following. 
By \Rem{result.for.final} the result is proven in the case 
$X$ admits a final object. To prove the 
result in general, one uses the descent properties of the objects involved:
that of the deformation functor according to~\ref{def.descent}, and 
that of the Deligne groupoid by~\cite{ddg}. The requirement on the existence of
a finite dimensional hypercover is due to a similar requirement in the 
proof of~\cite{ddg}, Theorem 4.1. 

Here are the details.

Choose $A$ to be a fibrant cofibrant algebra such that $A^i=0$ for $i>0$.
Choose a  finite dimensional hypercover $\epsilon: V_{\bullet}\to X$.

Let $A_n=A|_{V_n}$. This is a fibrant algebra on $X/V_n$; it is not necessarily
cofibrant but it ``behaves as if it were cofibrant''.

We have $T_A=\Right\Gamma(\CT_A)$. Let $T_n=\CT_A(V_n)=\der(A_n,A_n)$.

According to~\cite{ddg} (more precisely, according to its simplicial
version~\Prop{sddg}), there is an equivalence
\begin{equation*}
\label{}
\cdel_{T_A}\to\holim\{n\mapsto\cdel_{T_n}\}.
\end{equation*}

For each $n$ we define canonically a functor 
\begin{equation*}
\label{rhon}
\rho_n:\cdel_{T_n}\to\cdef_{A_n}
\end{equation*}
which will turn out to be an equivalence.

It is convenient to interpret here $\cdef_{A_n}$ as the homotopy
fiber of the map
\begin{equation*}
\label{}
\CW\to\ol{\CW}
\end{equation*}
where  $\CW$ denotes the simplicial groupoid $\widehat{\CW_*(R)}$
obtained by Dwyer-Kan localization from the simplical category $\CW_*(R)$
of $R\otimes\CO|_{V_n}$-algebras, and $\ol{\CW}=\widehat{\CW_*(k)}$.

This means, according to~\ref{hfib.scat}, that 
\begin{equation*}
\label{}
\cdef_{A_n}(R)=*\times_{\ol{\CW}}\ol{\CW}^{\Delta^1}\times_{\ol{\CW}}\CW.
\end{equation*}

The functor $\rho_n$ perturbs the differential. Here is its description.
Let $(R,\fm)\in\dgar(k)$. An object $z$ of $\cdel_{T_n}(R)$ is a
Maurer-Cartan element of $\fm\otimes T_n$. 
Then \\
$\rho_n(z)=(R\otimes A_n,d+z)$ where $d$ denotes the differential of 
the algebra $A$. This defines $\rho_n$ on objects.

Let now $\gamma:z\to z'$ be an $m$-morphism in $\cdel_{T_n}$. This means
that $\gamma\in\exp(\Omega_m\otimes\fm\otimes T_n)$ and $\gamma(z)=z'$.

Thus, $\gamma$ induces a map $\gamma: \rho_n(z)\to\Omega_m\otimes\rho_n(z')$
which gives an $m$-morphism in $\CW_*(R)$ and, therefore, in 
$\CW=\widehat{\CW_*(R)}$.

Let us check now that $\rho_n$ is an equivalence. Choose a cofibrant 
resolution $\pi:A'\to A_n$. Let $T'=\der(A',A')$. According to
\Lem{derrestriction-OK}, the maps
\begin{equation*}
\label{}
T'\overset{\phi}{\to} \der(A', A_n)\overset{\psi}{\leftarrow} T_n
\end{equation*} 
induced by $\pi$
are weak equivalences and, moreover, $\phi$ is surjective. Define
a dg Lie algebra $\fg$ by the cartesian diagram

\begin{center}
$$
\begin{diagram}[labelstyle=\scriptstyle ]
\fg & \rTo^{\phi'} & T_n \\
\dTo^{\psi'} &     & \dTo^{\psi} \\
T'  & \rTo^{\phi} & \der(A',A_n)
\end{diagram}.
$$
\end{center}

The maps $\phi'$ and $\psi'$
are, therefore, quasi-isomorphisms of dg Lie algebras.

We have the following functors:

\begin{itemize}
\item{ $\Psi:\cdel_{\fg}(R)\to \cdel_{T'}(R)$ induced by $\psi'$;}
\item{ $\Phi:\cdel_{\fg}(R)\to \cdel_{T_n}(R)$ induced by $\phi'$;}
\item{ $\rho_n:\cdel_{T_n}(R)\to \cdef_{A_n}(R)$ defined in~(\ref{rhon}).} 
\end{itemize}

The functors $\Psi$ and $\Phi$ are equivalences.
Define an arrow  $r:\cdel_{T'}(R)\to\cdef_{A_n}(R)$ as follows.

Recall that $\cdef_{A_n}(R)$ is the fiber of 
the map
\begin{equation*}
\label{}
\ol{\CW}^{\Delta^1}\times_{\ol{\CW}}\CW\lra\ol{\CW}.
\end{equation*}

Each object $z\in\cdel_{T'}(R)$ gives rise to an algebra $(R\otimes A',d+z)$
whose reduction is $A'$.
This gives an object in  $\cdef_{A_n}(R)$ presented by the morphism 
$A'\to A_n$. The action of $r$ on morphisms is obvious.

Look at the diagram
\begin{center}
$$
\begin{diagram}[labelstyle=\scriptstyle ]
\cdel_{\fg}(R) & \rTo^{\Psi} & \cdel_{T'}(R) & = & \cdef_{A'}(R) \\
\dTo^{\Phi}    &      & \dTo^{r} \\
\cdel_{T_n}(R) & \rTo^{\rho_n} & \cdef_{A_n}(R) 
\end{diagram}
$$
\end{center}

The diagram is not commutative. However, there is a homotopy
connecting $r\Psi$ with  $\rho_n\Phi$ assigning to each object 
$z\in\cdel_{\fg}(R)$
a morphism defined by the quasi-isomorphism
$$ (R\otimes A',d+\pi(z))\lra(R\otimes A_n,d+\chi(z))$$
induced by $\pi:A'\to A_n$.

The map $r$ is an equivalence (of homotopy fibers at $A_n$ and at $A'$).
Therefore, $\rho_n$ is also an equivalence.

The collection of functors $\rho_n$ sums up to an equivalence
\begin{equation*}
\label{totrho}
\holim(\rho):\holim\{n\mapsto \cdel_{\CT_A(V_n)}\}\to
\holim\{n\mapsto\cdef_{A_n}\}.
\end{equation*}  
We have already mentioned that the left hand side is equivalent to
$\cdel_{T_A}$. The right-hand side is equivalent to $\cdef_A$ by the descent
property~\ref{def.descent}. 
\end{proof}

\section{Examples}
\label{sect:examples}

\subsection{Deformations of schemes}

Let $X$ be a scheme over a field $k$ of characteristic zero.
Denote by $X^{\aff}_{\Zar}$ the site of affine open subschemes of $X$
with the Zariski topology. The topoi corresponding to the sites
$X_{\Zar}$ and $X^{\aff}_{\Zar}$ being equivalent, we can use  
$X^{\aff}_{\Zar}$ to describe deformations of $X$. The structure
sheaf $\CO$ considered as an object in $C((X^{\aff}_{\Zar})^{\pre}_k)$,
is fibrant. Therefore, our general result~\Thm{result} is applicable under some
finiteness conditions. We get the following.

\subsubsection{}
\begin{cor}{schemes}Let $X$ be a scheme over a field $k$ of characteristic
zero. Suppose that $X$ admits a finite dimensional hypercover by affine
open subschemes. Then the functor of formal deformations of $X$, $\cdef_X$,
is equivalent to the simplicial Deligne groupoid $\cdel_T$ where $T$ is
the dg Lie algebra of global derivations of a cofibrant resolution of $\CO$.
\end{cor}
\begin{proof}Let $p:A\rTo\CO$ be a  cofibrant resolution of $\CO$.
Since $\CO$ is fibrant and $p$ is an acyclic fibration, $A$ is fibrant as well.
The conditions~\ref{main-assumptions} are fulfilled, so \Thm{result}
gives the claim in question.
\end{proof}

The same reasoning provides a similar description
of formal deformations of a quasi-coherent sheaf of algebras.

\subsubsection{}
\begin{cor}{qcoh}
Let $X$ be a scheme over a field $k$ of characteristic
zero. Suppose that $X$ admits a finite dimensional hypercover by affine
open subschemes. Let $A$ be a quasi-coherent sheaf of algebras over a linear
operad $P$. Then the functor of formal deformations of $P$-algebra $A$
is equivalent to the simplicial Deligne groupoid $\cdel_T$ where $T$ is
the dg Lie algebra of global derivations of a cofibrant resolution of $A$.
\end{cor}

Note that here there is no connection between the $P$-algebra structure
and the $\CO_X$-module structure on $A$: we use the fact that $A$ is
quasi-coherent in order to deduce that $\CH^i(A)=0$ for positive $i$.

\subsection{Obstruction theory}
\label{obstruction.theory}
The tangent Lie algebra $T$ in Corollary~\ref{qcoh} is not easy to determine.
The classical obstruction theory task, the determination of the cohomology of 
$T$ is a much easier problem.

Assume we are dealing with deformations of associative algebras.

Let $A$ be a quasicoherent sheaf of associative $\CO_X$-algebras. 
We wish to describe deformations of $A$ as a $k$-algebra.

Consider $A$ as a presheaf on the site $X^{\aff}_{\Zar}$ of affine open 
subschemes of $X$ and let $\phi:P\rTo A$ be a cofibrant resolution. 
According to our definition, $\CT_A=\cpder(P,P)$ governs local deformations
of $A$. Let $C^*(P,P)$ be the Hochschild cochain complex for $P$.
Let $\pi:C^*(P,P)\rTo P$ be the obvious projection to $P=C^0(P,P)$.
The natural map
\begin{equation*}
\label{i_P}
i_P:\cpder(P,P)\to \Cone(\pi:C^*(P,P)\to P)
\end{equation*}
from the presheaf of derivations of $P$ to the (shifted and truncated) 
Hochshild cochains presheaf, is a weak equivalence.

According to Lemma~\ref{innerhoms}, $C^*(P,P)$ represents 
$\Right\cphom_{A\otimes A^{\op}}(A,A)$, the (local) Hochschild cohomology
of $A$. This gives the exact sequence
\begin{equation*}
\label{seelunts}
\ldots\rTo HH^i(A,A)\rTo H^i(X,A)\rTo H^i(T)\rTo HH^{i+1}(A,A)\rTo\ldots
\end{equation*}
where the global Hochschild cohomology $HH^i(A,A)$ is defined as
$$ H^i(\Right\Hom_{A\otimes A^{\op}}(A,A))$$
(compare to Lunts' \cite{lunts}, Cor. 5.4).

\subsection{Standard complex}

It seems too naive to expect that the standard complex 
$$\Cone(C^*(A,A)\rTo A))$$ 
represent the tangent Lie algebra $\CT_A$. 
However, in a very special case of associative deformations of the
structure sheaf of a smooth variety, the version of the standard complex
based on cochains which are differential operators in each argument,
gives a correct result.

In a more detail, let $A$ be the structure sheaf of a smooth algebraic 
variety $X$ over $k$, and let $C^*_{\DO}(A,A)$ denote the subcomplex of 
$C^*(A,A)$ consisting of cochains given by differential operators
in each argument.

Let us compare the Hochschild cochains $C^*(P,P)$ and $C^*_{\DO}(A,A)$.
One has the following commutative diagram
\begin{equation*}
\label{}
\begin{diagram}
& & T & & \\
& \ldTo && \rdTo^{\phi''} &\\
C^*(P,P) & & & & C^*_{\DO}(A,A)\\
& \rdTo^{\phi'} && \ldTo^{\psi} & \\
&& C^*(P,A) &&
\end{diagram}
\end{equation*} 
where $T$ is defined to make the diagram cartesian.
The map $\phi'$ is an acyclic fibration since $P^{\otimes n}$ are cofibrant
complexes of presheaves and $\phi$ is an acyclic fibration. Therefore,
$\phi''$ is an acyclic fibration of dg Lie algebras. By 
a version of Hochschild-Kostant-Rosenberg theorem proven in~\cite{yek},
$\psi$ is a weak equivalence of complexes of presheaves. This gives a weak
equivalence of Lie dg algebras $C^*(P,P)[1]$ and $C^*_{\DO}(A,A)[1]$.

\subsection{Equivariant deformations}
\label{equivariant.deformations}

Let $G$ be a group and let $\fg$ be a dg Lie algebra governing 
formal deformations of some object $A$ (we are informal at the moment).
If $G$ acts on $A$, one should expect a $G$-action to be induced on 
$\fg$. One can expect that the equivariant deformations of $A$ 
are governed by the dg Lie algebra $\Right\Gamma^G(\fg)$
whose $i$-th cohomology is $H^i(G,\fg)$. 

We are able to prove this when $G$ is a formal group, ``an object'' 
means ``a sheaf of operad algebras'', under the restrictions
of \Thm{result}.

In ~\ref{fg-beg}--\ref{fg-end} we discuss the action of formal groups on sites.
In \ref{equi.def} we describe the equivariant deformation functor.
The formula for the equivariant tangent Lie algebra is deduced 
in~\ref{equi-end}.

\subsubsection{}
\label{fg-beg}

In this subsection {\em a formal group} is a functor
$$ G:\art(k)\rTo \Groups$$
from the category of artinian local $k$-algebras to the category of groups,
commuting with the fiber products.

According to the formal Lie theory, the fiber at $1$ functor
\begin{equation*}
\label{fiber}
\ol{G}(R)=\Ker(G(R)\rTo G(k))
\end{equation*}
is uniquely defined by the corresponding Lie algebra $\fg$ (possibly, 
infinite-dimensional). The group $G(k)$ of $k$-points of $G$
acts on $\fg$ (adjoint action).

Thus, a formal group $G$ in our sense is described by a pair $(G(k),\fg)$
consisting of a discrete group $G(k)$, a Lie algebra $\fg$ and an action 
of $G(k)$ on $\fg$.

A representation of a formal group $G$ is a $k$-vector space $V$ together
with a collection of compatible operations 
$$ G(R)\rTo\ GL_R(R\otimes V).$$
This can be rephrased in terms of the corresponding pair $(G(k),\fg)$
as follows. To define a representation of the formal group corresponding
to $(G(k),\fg)$ on a $k$-vector space $V$, one has to define the action
of $G(k)$ and of $\fg$ on $V$ satisfying the compatibility condition
\begin{equation*}
\label{compatibility}
\gamma(x)(v)=\gamma(x(\gamma^{-1}v)), \quad \gamma\in G(k),\ x\in\fg,\ v\in V.
\end{equation*}

\subsubsection{}
Let $X$ be a site and $\CO$ be a sheaf of commutative $k$-algebras on $X$.
An action of a formal group $G$ on the ringed site $(X,\CO)$ is a collection
of compatible actions of groups $G(R)$ on $(X,R\otimes\CO)$ for $R\in\art(k)$.

If a formal group $G$ is described by a pair $({G(k)},\fg)$ as in~\ref{fg-beg},
its action on $(X,\CO)$ is given by an action of the discrete group ${G(k)}$
on $(X,\CO)$, action of $\fg$ by vector fields on $\CO(U)$ for each $U\in X$,
subject to the compatibility
\begin{equation*}
\label{compatibility.2}
\gamma(x)(\gamma(f))=\gamma(x(f)).
\end{equation*}

In the case $\CO=k$ (our assumption below) we will assume that the action
of $\fg$ on $k$ is trivial, so that the action of $G$ on $X$ is reduced to
the action of the discrete group $G(k)$.

\subsubsection{}
 Now we are able to define $G$-equivariant (pre)sheaves on $X$.
Suppose a formal group $G$ acts on $X$ as above. A $G$-(pre)sheaf
$M$ is given by a (pre)sheaf $M$ of $k$-vector spaces on $X$  together with 
a compatible collection of structures of $G(R)$-module on $R$-module 
$R\otimes M$.

If $G$ is presented by a pair $({G(k)},\fg)$ as above, a $G$-module
structure on $M$ amounts to a compatible collection
of maps $\gamma:M(U)\rTo M(^{\gamma}U)\quad (\gamma\in{G(k)})$, 
a collection of actions
$\fg\rTo\End(M(U))$ satisfying the condition
\begin{equation*}
\label{compatibility.3}
\gamma(xm)=\gamma(x)\gamma(m).
\end{equation*}
This implies the following

\begin{prop}{newsite}
There is an explicitly defined ringed site $(X/G,\CU)$ such that 
the category of $G$-(pre)sheaves of $k$-modules on $X$ is equivalent to
that of $(X/G,\CU)$-modules.
\end{prop} 
\begin{proof}
The category
$X/G$ has the same objects as $X$; for $U,V\in X$ one has
\begin{equation*}
\Hom_{X/G}(U,V)=\{(\gamma,f)|\gamma\in {G(k)},\,f\in\Hom_X(U,^\gamma V)\}.
\end{equation*}
The composition of morphisms in $X/G$ is defined in a standard way,
so that the morphism defined by a pair $(\gamma,f)$ is denoted as $\gamma f$.
One has the identity
$$ f\gamma=\gamma\,^\gamma f.$$
The topology on $X/G$ is generated by that on $X$.

We define the sheaf of rings $\CU$ as the sheafification of a presheaf
$\CU'$ defined below.
Let $U\in X$. We set $\CU'(U)$ to be the enveloping algebra of $\fg$.

If $f:U\rTo V$ is a map in $X$, the corresponding  map 
$$\CU'(V)\rTo\CU'(U)$$
is identity.
For $\gamma:^\gamma U\rTo U$ the corresponding map
$$ \CU'(U)\rTo \CU'(^\gamma U)$$
is induced by the automorphism of the Lie algebra $\fg$
sending $x\in\fg$ to $\gamma(x)$.

A straightforward check shows the pair $(X/G,\CU)$ satisfies the
required property.
\end{proof}
  
\subsubsection{}

A special case $X=*$ of the above construction gives the ringed site
$BG=(*/G,\CU)$. As a category, this is the classifying groupoid of ${G(k)}$;
the sheaf of rings is defined by the enveloping algebra $U\fg$ endowed with
the adjoint ${G(k)}$-action. Sheaves on this site are $G$-modules.

\subsubsection{} 
\label{fg-end}
If $M$ is a presheaf of complexes on $X/G$, we denote by 
$M^{\#}$ the presheaf on $X$ obtained by forgetting the $G$-structure. If $M$
is fibrant, $M^{\#}$ is fibrant as well; the global sections
$\Right\Gamma(X,M^{\#})$ admit a natural $G$-structure, that is define
a sheaf on $BG$ (this is just the higher direct image of the
morphism $X/G\rTo BG$ induced by the morphism $X\rTo *$). We denote this
sheaf $\Right\Gamma(X,M)$. One has 
\begin{equation*}
\label{Gglobalsec}
\Right\Gamma(X/G,M)=\Right\Gamma^G(\Right\Gamma(X,M)).
\end{equation*}

\subsubsection{Equivariant deformation functor}
\label{equi.def}

Let $X$ be a site and let a formal group $G$ act on $X$ (through $G(k)$).
Let $\CO$ be a $G$-equivariant operad. 
Let $\widehat\CW^\fl_{\she, G}(R),\ R\in\dgar(k),$
be the weak groupoid of $R$-flat sheaves of equivariant $R\otimes\CO$-algebras.
Let  $A$ be an $\CO$-algebra in $C(X^{\she},k)$ endowed with a $G$-action.

Similarly to \ref{def-functor}, we define the equivariant deformation 
functor
$$ \cdef_{A,G}:\dgar(k)\rTo\sGrp$$
as the homotopy fiber at $A$ of the functor 
$$ \widehat\CW^\fl_{\she, G}(R)\rTo \widehat\CW^\fl_{\she, G}(k)$$
induced by the projection $R\to k$.

\subsubsection{... and its tangent Lie algebra}
\label{equi.result}

We assume that the conditions of~\Thm{result} are satisfied for $\CO$-algebra
$A$. In particular, deformations of $A$ are governed by the global tangent
Lie algebra $T_A$ which can be calculated using a cofibrant resolution
$P$ of $A$ in the category $\Alg(\CO,X^{\pre}_k)$ by the formula
$$ T_A=\Right\Gamma(X,\cpder(P,P)).$$
We wish to express the functor of equivariant deformations through $T_A$.

This can be done as follows.
Consider the ringed site $(X/G,\CU)$ described in~\ref{newsite}.

Sheaves (resp., presheaves) on $(X/G,\CU)$ are precisely equivariant sheaves
(resp., presheaves) on $(X,k)$. 

We define a new operad $\CO\#G$ in $C((X/G)^{\pre}_k)$ as the one 
governing equivariant $\CO$-algebras on $X$. It is explicitly given by the 
formula
\begin{equation*}
\label{groupoperad}
\CO\#G(n)=\CO(n)\otimes\CU^{\otimes n},
\end{equation*}
with the operations uniquely defined by the $G$-action on $\CO$
\begin{equation*}
\label{comp.groupoperad}
\CU\otimes\CO(n)\rTo^{\Delta^n\otimes\id}\CU^{n+1}\otimes\CO(n)
\rTo\CU\otimes\CO(n)\otimes\CU^{\otimes n}\rTo\CO(n)\otimes\CU^{\otimes n},
\end{equation*}
where the second arrow swaps the arguments and the third one is defined
by the $G$-action on $\CO(n)$.

This is a generalization of the twisted group ring construction.

\subsubsection{}
\begin{lem}{XvsXG}
1. The forgetful functor $\#:(X/G)^{\pre}_{\CU}\to X^{\pre}_k$ preserves
weak equivalences, fibrations and cofibrations.

2. The same is true for the forgetful functor
$$\#:\Alg(\CO\#G,(X/G)^{\pre})\rTo \Alg(\CO,X^{\pre}).$$
\end{lem}
\begin{proof}
The proofs of both claims are identical. Since sheafification commutes 
with $\#$, weak equivalences 
are preserved. Cofibrations are preserved since joining a section over $U\in X$
corresponds, after application of $\#$, to joining sections
corresponding to a chosen basis of $U\fg$ over all 
$\gamma(U)$, for $\gamma\in G$. 
Finally, since any hypercover in $X/G$ is isomorphic to a hypercover in $X$, 
fibrations are also preserved.
\end{proof}

\subsubsection{}
\label{equi-end}
Thus, $A$ can be considered as a sheaf of
algebras on $X/G$. By \Lem{XvsXG} the condition $\CH^i(A)=0$ for $i>0$
is valid also when $A$ is considered as a sheaf on $X/G$. 
Therefore, \Thm{result} is applicable to $A$ in the equivariant setting.
Let us calculate the equivariant global tangent Lie algebra. Choose
an equivariant fibrant cofibrant representative $P$ of $A$. We can calculate
both equivariant and non-equivariant tangent Lie algebras using $P$.

Thus, the equivariant local tangent Lie algebra is 
\begin{equation*}
\label{Teq}
\CT_{A,G}=\cpder_{X/G}(P,P)
\end{equation*}
and the one in $X$ is
\begin{equation*}
\label{Tneq}
\CT_A=\cpder_X(P^{\#},P^{\#}).
\end{equation*}

Note that $\CT_A=\CT_{A,G}^{\#}$.
Since $\CT_{A,G}$ is fibrant by~\ref{der-fib}, one has
 
\begin{equation*}
\label{R=RR}
\Right\Gamma(X/G, \CT_{A,G})=\Right\Gamma^{G}\circ
\Right\Gamma(X, \CT_{A,G}).
\end{equation*}

\appendix\section{Simplicial categories and all that}
\label{App:scat}
%
%
%
%
%
%

\subsection{Simplicial categories}
\label{ss:cmc-scat}
Throughout the paper simplicial category means a simplicial object
in the category $\Cat$ of small categories having a discrete simplicial
set of objects. The category of simplicial categories is denoted by 
$\sCat$.

\subsubsection{}
The functor
$$\pi_0:\sCat\to\Cat$$
is defined by the formulas 
$$
\begin{array}{lll}
  \Ob \pi_0(X)&=&\Ob X\\
  \Hom_{\pi_0(X)}(x,y)&=&\pi_0(\Hom_X(x,y)).
\end{array}
$$
\subsubsection{}
Any simplicial category can be considered as a simplicial object
in $\Cat$. Applying the nerve functor
$$\CN:\Cat\to\simpl$$
layer by layer and taking the diagonal, we get a simplicial set
called {\em simplicial nerve} (or just nerve) $\CN(\CC)$ of a simplicial
category $\CC$.

\subsubsection{Model structure}
In this paper we use a model category structure
on $\sCat$ defined in~\cite{uc}.

\subsubsection{}
\begin{defn}{scat-we}
A map $f:\CC\to\CD$ in $\sCat$ is called a {\em weak equivalence}
if the map $\pi_0(f):\pi_0(\CC)\to\pi_0(\CD)$ induces a weak homotopy
equivalence of the nerves and for each $x,y\in\Ob(\CC)$ the map
$$ \shom(x,y)\to\shom(f(x),f(y))$$
of the simplicial Hom-sets is a weak equivalence. 
\end{defn}

Sometimes the following notion of {\em strong equivalence} is useful.

\subsubsection{}
\begin{defn}{scat-se}
A weak equivalence $f:\CC\to\CD$ is called {\em strong equivalence}
if the functor 
$\pi_0(f):\pi_0(\CC)\to\pi_0(\CD)$ is an equivalence of categories.
\end{defn} 

\subsubsection{}
Cofibrations in $\sCat$ are generated by the following maps:

(cof-1) $\emptyset\to *$, the functor from the empty simplicial category to a
one-point category.

(cof-2) For each  cofibration $K\to L$ in $\simpl$ the induced  map from
$K_{01}$ to $L_{01}$. Here $K_{01}$ denotes the simplicial category having 
two objects $0$ and $1$ with the only nontrivial maps $K=\shom(x,y)$.

\subsubsection{}
\begin{thm}{scat-cmc-uc}(\cite{uc})
The collections of cofibrations and of weak equivalences
define a CMC structure on $\sCat$.
\end{thm}

The maps in $\sCat$ satisfying RLP with respect to  acyclic cofibrations
(=cofibrations + weak equivalences) will be called fibrations.

Recall for the sake of completeness the explicit definition of fibration.

\subsubsection{}
\begin{defn}{scat-fib}A map $f:\CC\to\CD$ in $\sCat$ is called a fibration
if it satisfies the following properties

(1) the right lifting property (RLP) with respect to the maps
$$ \partial^{0,1}:\Delta^0\to\Delta^1$$
from the terminal category $\Delta^0=*$ to the one-arrow category $\Delta^1$.
 
(2) For all $x,x'\in\Ob\CC$ the map $f:\shom(x,x')\to\shom(fx,fx')$
is a Kan fibration. 
\end{defn}

\subsubsection{Simplicial structure and $\Tot$}
The category $\sCat$ admits a structure of a simplicial model category.

Let $S$ be a simplicial set, $\CC$ be a category
and let $\CN\CC$ be the nerve of $\CC$. The simplicial set $\shom(S,\CN\CC)$
is the nerve of a category which will be denoted by $\CC^S$.

Let $X=\{X_n\}\in\sCat,\quad S\in\simpl.$ The collection
$$ n\mapsto X_n^S$$
forms a simplicial object in $\Cat$. We define $X^S$ to be the
simplicial category given by the formulas
\begin{eqnarray}
  \label{eq:shom1}
  \Ob X^S=\Ob X_0^S\\
  \label{eq:shom2}
  \Hom_{X^S}(x,y)_n=\Hom_{X_n^S}(x_n,y_n)
\end{eqnarray}
where $x,y$ are objects of $X_0^S$ and $x_n,y_n$ are their degeneracies
in $X_n^S$.

Given a cosimplicial object $X^{\bullet}$ in $\sCat$, we can now define
$\Tot(X^{\bullet})$ by the usual formula
\begin{equation}
\label{tot-scat}
\Tot(X^{\bullet})=\invlim_{p\to q}(X^q)^{\Delta^p}.
\end{equation} 

Note that our simiplicial structure (and functor $\Tot$) on $\sCat$
is not standard. In~\cite{dhk} another defintion is given. The definitions
coincide for simplicial groupoids. However, the definition of~\cite{dhk}
does not provide $\sCat$ with a simplicial model structure, and 
this is the main reason we use the model structure described above.

\subsubsection{Homotopy limits}
\label{holim.scat}
The functor $\Tot$ defined in~(\ref{tot-scat}) prescribes a definition
of the homotopy limit as in, say, \cite{dhk}, Chapter XIV. If 
$F:I\to\sCat$ is a functor, $\holim F$ is defined as $\Tot(\widetilde{F})$
where the cosimplicial object  $\widetilde{F}$ in $\sCat$ is defined
by the standard formula
\begin{equation*}
\widetilde{F}^n=\prod_{i_0\to\ldots\to i_n\in\CN_n(I)}F(i_n).
\end{equation*} 

In the case all $F(i)$ are fibrant in $\sCat$, the homotopy limit 
$\holim\ F$ represents the right derived functor 
$\Right\lim:\Ho(\sCat^I)\to\Ho(\sCat)$.

In this case the following holds.


\subsubsection{}
\begin{prop}{holim-vs-N}(see~\cite{uc})
Let $F: I\to\sCat$ be a functor such that $F(i)$ are fibrant for all $i\in I$.
Then the natural functor 
$$ \CN(\holim F(i))\rTo \holim(\CN(F(i)))$$
is a weak equivalence.
\end{prop}

\subsubsection{Homotopy fibers}
\label{hfib.scat}
In this paper we are particularly interested in homotopy fibers. Let 
$f:\CC\to\CD$ be functor in $\sCat$ and let $d\in\CD$.  We suppose also that
the simplicial categories $\CC$ and $\CD$ are fibrant (for instance,
simplicial groupoids), so that the homotopy fiber represents the right
derived functor of the usual fiber product. In this case homotopy
fiber of $f$ at $d$ can be represented by the fiber product

\begin{equation*}
*\times_{\CD}\CD^{\Delta^1}\times_{\CD}\CC
\end{equation*}
where the map $*\to\CD$ is given by $d\in\CD$ and the maps from 
$\CD^{\Delta^1}$ to $\CD$ are given by the ends of the segment $\Delta^1$.

This immediately follows from~\cite{uc}, Prop. A.4.3 claiming, in particular,
that the map $\CD^{\Delta^1}\to \CD\times\CD$ is a fibration for fibrant $\CD$.

\subsection{Weak groupoids}
\label{weak-grp}

\subsubsection{}
\begin{defn}{hgr}A simplicial category $\CG$ is called a {\em 
weak groupoid} if $\pi_0(\CG)$ is a groupoid.
\end{defn}

The following fact justifies the above definition.
\subsubsection{}
\begin{prop}{wg-eq-gr}
A simplicial category $\CC$ is a weak groupoid if and only if it is
strongly equivalent to a simplicial groupoid. 
\end{prop}
\begin{proof}
The ``if'' part is obvious. According to~\cite{dk2}, if $\CC$ is a weak
groupoid, the map $\CC\to L^H(\CC,\CC)$ from $\CC$ to the hammock localization
of $\CC$ is a strong equivalence.
\end{proof}

\subsubsection{}
\begin{cor}{eq-vs-nerves}
A map $f:\CV\to\CW$ of weak groupoids is a weak equivalence
iff its nerve $\CN(f):\CN(\CV)\to\CN(\CW)$ is a weak homotopy equivalence. 
\end{cor}
\begin{proof}
\Prop{wg-eq-gr} reduces the claim to simplicial groupoids. In this case the
result is proven, for instance, in~\cite{dha}, 6.2.2, 6.2.3.
\end{proof}

\subsubsection{}
\begin{defn}{WGC}
Let $\CW$ be a simplicial category. Its weak groupoid completion
$\widehat{\CW}$ is the hammock localizaton $L^H(\CW,\CW)$, see~\cite{dk2}.
\end{defn}

One has a canonical map $\CW\to\widehat{\CW}$ from a simplicial
category to its weak groupoid completion. 

Since simplicial localizations preserve the homotopy type of the nerve,
see~\cite{dk}, 4.3, we obtain immediately the following

\subsubsection{}
\begin{cor}{eq-lcc-ldd}
A map $f:\CC\to\CD$ of simplicial categories induces a weak equivalence
of the weak groupoid completions $\widehat{\CC}\to\widehat{\CD}$ iff the nerve
$\CN(f):\CN(\CC)\to\CN(\CD)$ is a weak homotopy equivalence.
\end{cor}

\subsection{Weak groupoid of a model category}

Let $\CC$ be a closed model category and let 
$\CW$ be the subcategory of weak equivalences in $\CC$.

One can assign to $\CC$ a few different weak groupoids.
These are weak groupoid completions of the categories $\CW$, $\CW^{\c}$
(the objects are cofibrant objects of $\CC$, the morphisms are weak
equivalences), and also $\CW^{\f}$ and $\CW^{\cf}$. We denote them
$\widehat{\CW}$, $\widehat{\CW^{\c}}$, $\widehat{\CW^{\f}}$ and 
$\widehat{\CW^{\cf}}$.  

In the case $\CC$ admits a simplicial structure so that the axiom 
(SM7) of~\cite{ha} is satisfied, one defines another weak
groupoid $\CW^{\cf}_*$ whose objects are cofibrant fibrant objects of $\CC$
and $n$-morphisms are the ones lying in the components of weak equivalences.
Similarly the weak groupoid $\widehat{W_*}$ is defined.

In \Prop{eq-weak-gr} below we show that all these weak groupoids are strongly 
equivalent.

\subsubsection{Contractibility of resolutions}

Homological algebra starts with an observation that 
resolutions are usually unique up to a homotopy which is itself unique up to
homotopy. In this section we prove a generalization of this fact:
the category of resolutions has a contractible nerve. This result will
be used in the proof of equivalences~\ref{eq-weak-gr} below. 

Let $\CC$ be a closed model category, $M\in\CC$. Let $\CC^{\c}_M$ 
denote the full subcategory of the category $\CC/M$ whose objects are the 
weak equivalences $f:P\to M$ with cofibrant $P$. Similarly, $\CC^{\cf}_M$
consists of quasi-isomorphisms $P\to M$ with $P$ cofibrant and fibrant.

\subsubsection{}
\begin{thm}{res-cont}1. The nerve of the category $\CC^{\c}_M$ is contractible.

2. If $M$ is fibrant, the nerve of $\CC^{\cf}_M$ is contractible.
\end{thm}
\begin{proof}
{\em Step A}. First of all, we check that the nerve of $\CC^{\c}_M$ is simply 
connected. Suppose $f:P\to M$ and $g:Q\to M$ are two objects of $\CC^{\c}_M$.
Present the map $f\coprod g:P\coprod Q\to M$ as a composition of a cofibration
$\iota$ and an acyclic fibration $\pi$
\begin{equation*}
\label{}
P\coprod Q\overset{\iota}{\to} R\overset{\pi}{\to} M.
\end{equation*}
The map $\pi:R\to M$ presents an object of $\CC^{\c}_M$. This proves 
that the nerve of $\CC^{\c}_M$ is connected. Note that the same construction
proves that the nerve of $\CC^{\cf}_M$ is connected if $M$ is fibrant.

To prove that the nerves in question are simply connected, one can pass
to groupoid completions and calculate the automorphism group of any
object of the obtained groupoid. A standard reasoning shows that an
acyclic fibration $f:P\to M$ with cofibrant $P$ has no nontrivial
automorphisms in the groupoid completion. 

{\em Step B}. Choose an acyclic fibration $f:P\to M$ with cofibrant $P$.
This defines a functor
\begin{equation*}
\label{fibre}
F:\CC^{\c}_P\to\CC^{\c}_M
\end{equation*}
which carries a weak equivalence $g:Q\to P$ to $fg:Q\to M$.
Fix an object $g:R\to M$ of $\CC^{\c}(M)$. The fibre category $F/g$
can be easily identified with the category $\CC^{\c}_{P\times_MR}$.
According to Step A, all fibres of $F$ are simply connected.
Note that the nerve of $\CC^{\c}(P)$ is contractible since the category
admits a final object.

{\em Step C}. We have to check that the reduced homology of the nerve
of $\CC^{\c}_M$ vanishes. The construction of Step B allows to prove this by 
induction. In fact, let $\widetilde{H}_i(\CC^{\c}_M)=0$ for all $i<n$
and for all objects $M$. Then \Prop{ThmC} below shows that 
$\widetilde{H}_n(\CC^{\c}_M)=0.$ This proves the theorem modulo \Prop{ThmC}
below.
\end{proof}

The following result is very much in the spirit of~\cite{hakt1}, Theorems A 
and B.

\subsubsection{}
\begin{prop}{ThmC}Let $F:\CC\to\CD$ be a functor. Suppose that

(a) the nerve of $\CC$ is contractible.

(b) $\widetilde{H}_i(F/d)=0$ for all $i<n,\ d\in \CD$.

Then $\widetilde{H}_n(\CD)=0$.
\end{prop}
\begin{proof}
Consider the bisimplicial set $T_{\bullet\bullet}$ (used by Quillen
in the proof of Theorem A, cf.~\cite{hakt1}, p.~95) defined by the formula

\begin{equation*}
\label{}
T_{pq}=\{c_q\to\ldots c_0;\, F(c_0)\to d_0\to\ldots\to d_p\}, 
\end{equation*}
with the obvious faces and degeneracy maps.

The diagonal of this bisimplicial set, $\diag\ T$, is homotopy 
equivalent to the nerve of $\CC$ (see~\cite{hakt1}, p.~95). Therefore, it has
trivial homology.

Recall that the homology of a simplicial set $X$ can be calculated as follows.
First, one considers $\mathbb{Z}X\in\Dop\Ab$. Then one applies
Dold-Puppe equivalence of categories
\begin{equation*}
\Norm:\Dop\Ab\to C^{\leq 0}(\mathbb{Z}).
\end{equation*}
Finally, one has $H_i(X)=H^{-i}(\Norm(\mathbb{Z}X))$.

If $Y\in (\Dop)^2\Ab$, we denote $\Norm^2(Y)$ the bicomplex obtained from $Y$
by normalization in both directions. 

\begin{lem}{}
$$\Norm(\diag(Y))=\Tot(\Norm^2(Y)).$$
\end{lem}
The lemma is similar to Quillen's lemma at p.~94, ~\cite{hakt1}. One checks
it for representable $Y=\mathbb{Z}h^{pq}$ where
$h^{pq}_{rs}=\Delta^p_r\times\Delta^q_s$, and then checks that both sides
of the equality commute with the direct limits.

Now consider the map
\begin{equation*}
\label{}
T_{pq}=\coprod_{d_0\to\ldots\to d_p}\Nerve(F/d_0)\to
\coprod_{d_0\to\ldots\to d_p}\pt=\Nerve_p(\CD).
\end{equation*}

Let $Z$ be the bicomplex corresponding to the bisimplicial abelian group 
$\mathbb{Z}T_{\bullet\bullet}$. Denote by $H^{\vert}$ and $H^{\hor}$
the homology with respect to the vertical (of degree $(0,1)$) and the 
horizontal (of degree $(1,0)$) differential. By the assumptions of the
proposition, one has $H^{\vert}_q(Z)=0$ for $q=1,\ldots,n-1$
and $H^{\vert}_0(Z)=\Norm_p(\Nerve(\CD))$. Look at the spectral sequence
\begin{equation*}
\label{}
E^2_{pq}=H^{\hor}_pH^{\vert}_q(Z)\Rightarrow H_n(\Tot(Z)).
\end{equation*}
According to the above lemma, the spectral sequence converges to zero. 
Our calculation shows that $E^2_{pq}=0$ for $q=1,\ldots,\, n-1$ and 
$E^2_{p0}=H_p(\Nerve(\CD))$. This implies \Prop{ThmC} since the map
$H_n(\Nerve(\CD))=E^2_{n0}\to E^{\infty}_n$ should be injective.
\end{proof}
\subsubsection{}
\begin{prop}{eq-weak-gr}
(i) The weak  groupoids $\widehat{\CW}$, $\widehat{\CW^{\c}}$, 
$\widehat{\CW^{\f}}$ and $\widehat{\CW^{\cf}}$ are equivalent.

(ii) Suppose $\CC$ admits a simplicial structure and suppose that (SM7)
and half the axiom (SM0) of ~\cite{ha},  (existence
of simplicial cylinder or path spaces), is fulfilled. Then the weak
groupoids  $\CW^{\cf}_*$ and $\widehat{W_*}$ are also equivalent to 
the above.
\end{prop}
\begin{proof}
(i) According to~\ref{eq-vs-nerves} it is enough to prove the nerves of
the four categories are weakly homotopically equivalent. This follows
from~\Thm{res-cont} by Quillen's Theorem A, see~\cite{hakt1}.

(ii) Since $\CW^{\cf}_*$ is a weak groupoid, it is equivalent to its
weak groupoid completion. In order to prove that the latter is equivalent
to the weak groupoid completion of $\CW^{\cf}$, one has to compare the
(simplicial) nerves of $\CW^{\cf}_*$ and of $\CW^{\cf}=\CW^{\cf}_0$.
This will immediately follow once we prove that the multiple
degeneracy $s:\CW^{\cf}_0\to\CW^{\cf}_n$ induces a weak equivalence 
$\CN(s):\CN(\CW^{\cf}_0)\to\CN(\CW^{\cf}_n)$.

Let us suppose that the simplicial path functor exists. Then
$\Hom_n(x,y)=\Hom(x,y^{\Delta^n})$. Therefore, the fibre $s/y$
is equivalent to the category of fibrant cofibrant resolutions
of the object $y^{\Delta^n}$. This object being fibrant, \Thm{res-cont}
asserts that $s/y$ is contractible. Once more Quillen's Theorem A
accomplishes the proof. The case when the simplicial cylinder functor exists,
as well as the weak groupoid $\widehat{W_*}$ are treated similarly.
\end{proof}

\subsection{Simplicial Deligne groupoid}
\label{deligne}
In this paper weak groupoids appear as values of a formal 
deformation functor on artinian algebras. One looks for a presentation 
of such a functor with a dg Lie algebra. We recall below three functors
assigned to a dg Lie algebra, of which the last one is used in this paper.

\subsubsection{Deligne groupoid \em{(see~\cite{gm})}}
Let  $\fg$ be a Lie dg algebra and $(R,\fm)\in\dgar(k)$. Deligne groupoid
$\Del_{\fg}(R)$ has as objects the Maurer-Cartan elements of $\fm\otimes\fg$,
\begin{equation*}
\label{Ob-d-g}
\Ob\Del_{\fg}(R)=\MC(\fm\otimes\fg):=
\{z\in (\fm\otimes\fg)^1|dz+\frac{1}{2}[z,z]=0\}.
\end{equation*}
The group $\exp(\fm\otimes\fg)^0$ acts in a natural way on the 
set~(\ref{Ob-d-g}) and one defines
\begin{equation*}
\label{Mor-d-g}
\Hom_{\Del_{\fg}(R)}(z,z')=\{\gamma\in\exp(\fm\otimes\fg)^0|z'=\gamma(z)\}.
\end{equation*}

This definition is homotopy invariant if one requires $(\fm\otimes\fg)^i=0$ 
for $i<0$. 

\subsubsection{Nerve of a dg Lie algebra {\em (see~\cite{ddg,uc})}}
Let  $\fg$ and $R$ be as above. Define for $n\geq 0$ 
\begin{equation*}
\label{nerve-dgl}
\Sigma_{\fg}(R)_n=\MC(\Omega_n\otimes\fm\otimes\fg)
\end{equation*} 
where $\Omega_n$ is the algebra of polynomial differential forms on the
standard $n$-simplex. Then $\Sigma_{\fg}(R)$ is a Kan simplicial set,
its fundamental groupoid is canonically identified with $\Del_{\fg}(R)$
and the two are homotopically equivalent if $\fm\otimes\fg$ is 
non-negatively graded.

\subsubsection{Simplicial Deligne groupoid {\em (see~\cite{dha})}}
Here $\fg$, $R$ are as before. Define simplicial groupoid $\cdel_{\fg}(R)$
as follows. Its objects are the objects of $\Del_{\fg}(R)$. The collection
of $n$-morphisms from $z$ to $z'$ coincides with the collection of 
such morphisms in the Deligne groupoid corresponding to the Lie dg algebra 
$\Omega_n\otimes\fg$.

The following lemma connects between the different constructions.
\subsubsection{}
\begin{lem}{simp-del-nerve}
1. There is a natural weak equivalence of simplicial sets
$$\CN(\cdel_{\fg}(R))\isom \Sigma_{\fg}(R).$$

2. Deligne groupoid $\Del_{\fg}$ naturally identifies with the
fundamental groupoid $\Pi(\Sigma_{\fg})$ of the nerve $\Sigma_{\fg}$.
The composition 
$$\Sigma_{\fg}(R)\to\Pi(\Sigma_{\fg}(R))=\Del_{\fg}(R)$$
is a weak equivalence if $\fm\otimes\fg$ is non-negatively graded.
\end{lem}

\subsubsection{Descent}

All three functors mentioned above are defined by a nilpotent dg Lie algebra
$\fg_R:=\fm\otimes\fg$: one has $\Del_{\fg}(R)=\Del(\fg_R)$,
$\Sigma_{\fg}(R)=\Sigma(\fg_R)$ and $\cdel_{\fg}(R)=\cdel(\fg_R)$
in an obvious notation. In what follows we will fix once an forever
the commutative dg algebra $R$ and we will erase the subscript $R$ from the
notation.

Let $\fg^{\bullet}$ be a nilpotent cosimplicial dg Lie algebra.
A natural morphism of  simplicial sets
\begin{equation*}
\label{map-ddg}
\Sigma(\Tot(\fg^{\bullet}))\to\Tot(\Sigma(\fg^{\bullet}))
\end{equation*}
can be easily constructed. The main result of~\cite{ddg} claims
that this map is a homotopy equivalence provided $\fg^{\bullet}$ is
finitely dimensional in the cosimplicial direction.

In the main body of the paper we need a similar result in the context of
simplicial Deligne groupoids. Let us show it easily follows from
the result of~\cite{ddg}.

Let us construct a map of simplicial groupoids
\begin{equation}
\label{map-ddg-sg}
\cdel(\Tot(\fg^{\bullet}))\to\Tot(\cdel(\fg^{\bullet})).
\end{equation}

On the level of objects the map is constructed as follows.
An object of the left-hand side is an element of
$$ \invlim_{p\to q}\MC(\Omega_p\otimes\fg^q).$$
An element of $\MC(\Omega_p\otimes\fg^q)=\Sigma_p(\fg^q)$ 
gives rise to a sequence of $p$ 1-simplices in $\Sigma(\fg^q)$; passing
to the fundamental groupoid we get a $p$-simplex in the groupoid $\Del(\fg^q)$.
Since an object of the right-hand side of~(\ref{map-ddg-sg}) 
is an element of 
$$ \invlim_{p\to q}\Ob\ \Del(\fg^q)^{\Delta^p},$$
the morphism~(\ref{map-ddg-sg}) is defined on the level of objects.

Fix $z,z'\in\MC(\Tot(\fg^{\bullet}))$.
A $n$-map from $z$ to $z'$ on the left-hand side of~(\ref{map-ddg-sg})
is given by an element $\gamma\in\exp(\Omega_n\otimes\Tot(\fg^{\bullet}))^0$
satisfying the equation $z'=\gamma(z)$. 

The composition
\begin{eqnarray*}
  \label{eq:composition}
\exp(\Omega_n\otimes\Tot(\fg^{\bullet}))^0\to
\exp(\Tot(\Omega_n\otimes\fg^{\bullet}))= 
\invlim_{p\to q}\exp(\Omega_p\otimes\Omega_n\otimes\fg^q)^0\to\\
\to\invlim_{p\to q}\CN_p(\Del(\Omega_n\otimes\fg^q))
\end{eqnarray*}
defines the map (\ref{map-ddg-sg}) for the $n$-morphisms.

\begin{prop}{sddg}
Suppose $\fg^{\bullet}$ is finite dimensional in the cosimplicial direction,
i.e. there exists $n$ such that the intersection of kernels
of all codegeneracies vanishes in degrees $>n$. Then the 
morphism~(\ref{map-ddg-sg}) defined above is an equivalence.
\end{prop}

To prove a map of simplicial groupoids is an equivalence, it is enough
to check it induces a homotopy equivalence of the nerves. Applying the 
nerve functor to the both sides of~(\ref{map-ddg-sg}), we get the 
morphism~(\ref{map-ddg}) which is an equivalence by~\cite{ddg}.
This proves the proposition.

\appendix\section{Simplicial presheaves}

\label{App:sim}

A model category structure similar to the one described in 
Section~\ref{sect:models} exists also on the category of simplicial 
presheaves. This model category structure differs from the one defined 
in~\cite{j}. More precisely, weak equivalences are the same; we have much 
less cofibrations and, consequently, much more fibrations.

Let $X$ be a site. Since the category $\simpl$ of simplicial sets
is cofibrantly generated, the category of simplicial presheaves
$\Dop(X^{\pre})$ admits a CMC structure in which cofibrations are
generated by cofibration in $\simpl$, see~\cite{dhk}, 9.6 and 6.10.

In this model structure a map $f:A\to B$ is a weak equivalence
(resp., a fibration) iff for each $U\in X$ $f(U):A(U)\to B(U)$ is a 
weak equivalence (resp., a fibration). The collection of cofibration is 
generated by gluing simplices along a boundary over some $U\in X$.

This is the model structure we have in mind in the case $X$ is endowed with
the coarse topology. In order to describe the model structure on 
$\Dop(X^{\pre})$ which remembers the topology of $X$, we redefine the
notion of weak equivalence as in~\cite{j}, leaving the notion of cofibration
unchanged. Recall the following definition due to Joyal~\cite{jo} and
Jardine, \cite{j}.
 
For $A\in\Dop X^{\pre}$ one defines $\pi_0(A)$ as the sheafification
of the presheaf $U\mapsto\pi_0(A(U))$. Similarly, for $n>0$ and 
$a\in A(U)_0$ one defines $\pi_n(A,a)$ as the sheafification of the presheaf
$V\mapsto\pi_n(A(V);a)$ defined on $X/U$.

\subsection{}
\begin{defn}{swe}
A map of simplicial presheaves $f:A\to B$ is a weak equivalence if
  \begin{itemize}
  \item $\pi_0(f):\pi_0(A)\to \pi_0(B)$ is an isomorphism;
  \item for each $a\in A(U)_0$ the map $\pi_n(f):\pi_n(A;a)\to\pi_n(B;f(a))$
is an isomorphism.
  \end{itemize}
\end{defn}

\subsection{}
\begin{thm}{simpl}The category $\Dop(X^{\pre})$ endowed with the classes
of weak equivalences described above and cofibrations as for the coarse
topology, is a closed model category.
\end{thm}

A nice feature of this model structure is the following description
of fibrations (see \Prop{s:fibrations} below).

\begin{prop}{}
A map $f:M\to N\in\Dop(X^{\pre})$ is a fibration if and only if the following
conditions are satisfied.
\begin{itemize}
\item{$f(U):M(U)\to N(U)$ is a Kan fibration for each $U\in X$.}
\item{For each hypercover $\epsilon:V_{\bullet}\to U$ the comutative
diagram
\begin{center}

\begin{diagram}
M(U)& \rTo & \Cech(V_{\bullet},M):=\Tot M(V_{\bullet})\\
\dTo&      &  \dTo                                    \\
N(U)& \rTo & \Cech(V_{\bullet},N)
\end{diagram} 
\end{center}
is homotopy cartesian.  }
\end{itemize} 
\end{prop}

\subsubsection{}
\begin{rem}{split-hc}
The Cech complex appearing in the description of fibrations, is not
homotopy invariant, even for pointwise fibrant presheaves. This means that
that fibrantness is not necessarily preserved under pointwise weak
equivalence of pointwise fibrant presheaves. 

A pointwise fibrant presheaf $F$ is called to satisfy 
{\em descent property} with respect to a hypercover $V_{\bullet}\to U$ if
$F(U)$ is homotopy equivalent to $\holim F(V_{\bullet})$.

In general, our fibrant presheaves do not satisfy the descent property.
However, suppose $V_{\bullet}$ is {\em split}, 
see~\cite{sga4}, \'Exp. Vbis, 5.1.1. This means that for each $n$
the map $D_n\to V_n$ from the subpresheaf of degenerate $n$-simplices
to the presheaf of all $n$-simplices, can be presented as a composition
$$ D_n\to D_n\sqcup N_n\isom V_n.$$
Then, for each
acyclic cofibration $S\to T$ of simplicial sets the induced map
of simplicial presheaves
$$ (S\times V_n)\coprod^{S\times D_n} (T\times D_n)\to T\times V_n$$
is an acyclic cofibration. This implies that if $F$ is a fibrant presheaf
and $V_{\bullet}$ is a split hypercover, then $F(V_{\bullet})$ is a fibrant
cosimplicial simplicial set in the sense of~\cite{bk}, X.4. This implies
that $\Tot F(V_{\bullet})$ is homotopy equivalent to $\holim F(V_{\bullet})$
and, therefore, $F$ satisfies the descent property with respect to 
$V_{\bullet}$. This implies that if any hypercover in $X$ is refined by
a split hypercover, our model structure coincides with the one named
$U\CC/S$ in~\cite{dhi}, Theorem 1.3.
\end{rem}

\subsubsection{}
To prove the theorem, we describe a collection of morphisms which are
simultaneously weak equivalences and cofibrations. These morphisms
are called {\em generating acyclic cofibrations}. \Thm{simpl} then follows
from \Lem{s:af-are-correct} below claiming that weak equivalences 
satisfying the RLP with respect to all generating acyclic cofibrations, 
are pointwise acyclic Kan fibrations.

\subsubsection{}The following version of~\cite{sga4}, V.7.3.2, plays a very
important role here.

Let $E$ be a topos, $M,N\in\Dop E$. We will use the notion of weak equivalence
in $\Dop E$ defined in~\ref{swe}.

\begin{prop}{asSGA4}
Let $f:M\to N$ be a morphism of simplicial objects in $E$. Suppose that 
\begin{itemize}
\item $f_p$ is an isomorphism for $p<n$.
\item $f_n$ is an epimorphism.
\item morphisms $M\to \cosk_n(M),\ N\to \cosk_n(N)$, are isomorphisms.
\end{itemize}
Then $f$ is a weak equivalence.
\end{prop}
\begin{proof}
The proof is essentially the same as in~\cite{sga4}, V.7.3.2.
Embed $E$ into a topos $E'$ having enough points. If $E=X^{\she}$,
$E'$ can be taken to be $X^{\pre}$. Let 
$$ a^*:E'\rlarrows E: a_*$$
be the corresponding inverse and direct image functors.
We construct $f':M'\to N'$ as follows. 
\begin{eqnarray*}
M'& = & a_*(M).\\
N'_i & = & a_*(N_i)\text{ for }i<n;\,  N'_n=a_*(f)(M'_n);\, N'=\cosk_n(N').\\
f' & = &a_*(f):M'\to N'.  
\end{eqnarray*}
Then $f=a^*(f')$.

The functor $a^*$ preserves weak equivalences.
In fact, weak equivalences can be described using finite limits and 
arbitrary colimits which are preserved by the inverse image functor.

This reduces the claim to the case $E$ has enough points.
In this case  $f$ is a weak equivalence iff for each point $x$ $f_x$ is a weak
equivalence. Inverse image functor preserves the properties listed in the 
proposition. Thus everything is reduced to the case $E=\Ens$ where the claim 
is well-known. 

\end{proof}

\subsubsection{}We denote by $\Delta_+$ the category obtained from $\Delta$
by attaching an initial object $\emptyset=[-1]$. In particular,
any hypercover $\epsilon:V_{\bullet}\to U$ defines an object  
$V_{\circ}\in\Delta_+^{\op}X^{\pre}$ with $V_{-1}=U$. In what follows
we will use the subscript (resp., the superscript) $\circ$ to denote the
augmented simplicial (resp., cosimplicial) object.

Let $A_{\bullet}^{\circ}\in\Delta_+(\Dop\Ens)$ and let 
$V_{\circ}\in\Delta_+^{\op}X^{\pre}$.

We define a simplicial presheaf $A_{\bullet}^{\circ}\otimes V_{\circ}$
by the formula
\begin{equation*}
\label{AotimesV}
(A\otimes V)_n=\dirlim_{p\to q\in\Delta_+}A^p_n\times V_q.
\end{equation*}

\subsubsection{}Let us make a few calculations.
 Let, for instance, $A^p_n=K_n$ for some $K\in\simpl$.
Then $A\otimes V=K\times V_{-1}$. 

\subsubsection{}
Another important example is
\begin{equation*}
\label{ktot}
A^p_n=K_n\times\Delta^p_n.
\end{equation*}

One can easily see that $(\Delta\otimes V_{\circ})=V_{\bullet}$.
Also, if $K\in\Dop\Ens$  then for any $ A^{\circ}_{\bullet}$ one has
\begin{equation*}
\label{timesotimes} 
(K\times A^{\circ}_{\bullet})\otimes V_{\circ}=
K\times (A^{\circ}_{\bullet}\otimes V_{\circ}).
\end{equation*}

Thus we finally get
\begin{equation*}
\label{KD}
(K\times\Delta)\otimes V_{\circ}=K\times V_{\bullet}.
\end{equation*}

An easy calculation shows that for a simplicial presheaf $M$ one has
\begin{equation*}
  \label{eq:ktot}
  \Hom(A\otimes V_{\circ},M)=\Hom(K,\Cech(V_{\bullet},M)).
\end{equation*}

\subsubsection{}Generating acyclic cofibrations for our model structure
consist of two collections. 

The first collection is numbered by $U\in X$ and
a pair of integers $(i,n)$ such that $0\leq i\leq n$. It consists of maps
\begin{equation*}
  \label{eq:ac1}
  \Lambda^n_i\times U\to\Delta^n\times U.
\end{equation*}
This collection defines the model category structure on the simplicial
presheaves on $X$ corresponding to the coarse topology.

The second collection is numbered by hypercovers $V_{\circ}$ of $X$ 
and integers $n\geq 0$. It is defined as follows.

Let $\Delta^+$ (don't confuse with $\Delta_+$!) denote the cosimplicial 
simplicial space with
\begin{eqnarray*}
  \label{eq:+}
  \Delta^{+n}=\Delta^{n+1}\\
  \delta^{+i}=\delta^i\\
  \sigma^{+i}=\sigma^i.
\end{eqnarray*}
The map $i:\Delta\to\Delta^+$ is defined by the ``last face'':
$$ i^n=\delta^{n+1}:\Delta^n\to\Delta^{n+1}.$$
  
Define $A(n)$ by the cocartesian diagram
\begin{equation}
\label{An}
\begin{diagram}
\partial\Delta^n\times\Delta & \rTo^i &  \partial\Delta^n\times\Delta^+ \\
\dTo                         &      &  \dTo                           \\
\Delta^n\times\Delta         & \rTo & A(n)
\end{diagram}
\end{equation}  

and put $B(n)=\Delta^n\times\Delta^+$. 
Then the natural map $A(n)\otimes V\to B(n)\otimes V$ is the basic acyclic 
cofibration corresponding to the pair $(V,n)$.

\subsubsection{}
The object $\Delta^+\in\Delta_+\Dop\Ens$ identifies with the pushout of 
the diagram
$$
\begin{diagram}
\Delta & \rTo^{\delta^0\times\id_{\Delta}} & \Delta^1\otimes\Delta \\
\dTo   &         & \dTo                  \\
*      & \rTo    & \Delta^+
\end{diagram}
$$
Under this identification the map $i:\Delta\to\Delta^+$ is induced by
$$\delta^1\times\id_{\Delta}:\Delta\to\Delta^1\times\Delta.$$

Therefore, $\Delta^+\otimes V_{\circ}$ can be calculated 
from the cocartesian diagram
$$
\begin{diagram}
V_{\bullet} & \rTo^{\delta^0\times\id_{V_{\bullet}}} & 
                               \Delta^1\times V_{\bullet} \\
\dTo        &                 & \dTo                       \\
V_{-1}      & \rTo^v            & \Delta^+\otimes V_{\circ}
\end{diagram}
$$

\begin{lem}{Cone}
The morphism $v: V^{-1}\to  \Delta^+\otimes V_{\circ}$ is a weak equivalence.
\end{lem}
\begin{proof}The simplicial set $\Delta^1$ being contractible,
the map $\delta^0\times\id_{V_{\bullet}}$ is a weak equivalence. It is also
pointwise injective. This means that it is a trivial cofibration in the 
sense of Jardine~\cite{j}. Then Proposition 2.2 of~\cite{j} asserts that
$v$ is a weak equivalence.
\end{proof}

\subsubsection{}
\begin{lem}{bac.are.ac}
Let $V_{\circ}$ be a hypercover and $n\in\mathbb{N}$. 
The map $A(n)\otimes V_{\circ}\to 
B(n)\otimes V_{\circ}$ of simplicial presheaves is a cofibration 
and a weak equivalence.
\end{lem}
\begin{proof}
The map is cofibration since $V_i$ are coproducts of representable 
presheaves. Let us check the acyclicity.

Tensoring diagram~(\ref{An}) by $V_{\circ}$ and using Proposition 2.2 
of~\cite{j} we get that $A(n)\otimes V_{\circ}$ is equivalent to $V_{\bullet}$.
The tensor product $B(n)\otimes V_{\circ}$ is obviously equivalent to $V_{-1}$.

By \Prop{asSGA4} the map $\epsilon:V_{\bullet}\to V_{-1}$ 
is a weak equivalence. This implies the lemma. 
\end{proof}

The following lemma is the analog of~\Lem{af}.
\subsubsection{}
\begin{lem}{s:af-are-correct}
Let $f:M\to N$ be a weak equivalence. Suppose $f$ satisfies
the RLP with respect to all generating acyclic cofibrations. Then
$f(U):M(U)\to N(U)$ is an acyclic Kan fibration for all $U\in X$.
\end{lem}
\begin{proof}
The map $f(U)$ is a fibration for any $U\in X$. Fix $U$ and fix a section 
of $N$ over $U$. Let $F$ be the fiber of $f$ at the chosen section. 
We wish to prove that $F(U)$ is acyclic. To prove this we have to check
that any commutative diagram
\begin{center}
$$
\renewcommand{\labelstyle}{\scriptscriptstyle}
  \begin{diagram}
  \partial\Delta^n & \rTo      & F(U)                    \\
  \dTo             & \ruDotsto &                            \\
  \Delta^n         &       &
  \end{diagram}
$$
\end{center}
can be completed with a dotted arrow. Below we denote $A(n)$ by $A$.

Since $f$ is a weak equivalence, there exists a covering $\epsilon:V_0\to U$ 
and a dotted arrow $A^0\to F(V_0)$ making the  diagram 

\begin{center}
$$
\renewcommand{\labelstyle}{\scriptscriptstyle}
  \begin{diagram}
  A^{-1}=\partial\Delta^n & \rTo      & F(U)& \rTo &F(V_0)     \\
  \dTo                    &           &     & \ruDotsto(4,2) &           \\
  A^0                     &           &     &      &  
  \end{diagram}
$$
\end{center}
commutative.

Suppose, by induction, a $k$-hypercover $\epsilon:V_{\leq k}\to U$
and a collection of compatible maps $A^i\to F(V_i),\ i\leq k$, has been 
constructed. This induces a map $\sk^k(A^{\bullet})^{k+1}\to F(V_k)$.
Since the sheafification of the homotopy groups $\pi_i(F)$ vanishes,
there exists a covering $V_{k+1}\to\cosk_n(V)_{n+1}$ and a map
$A^{k+1}\to F(V_{k+1})$ compatible with he above.

Therefore, a hypercover $\epsilon:V_{\bullet}\to U$ and a map 
$A(n)\otimes V\to F$ is constructed. According to hypothesis of the lemma,
this latter can be extended to a map $B(n)\otimes V\to F$ which includes
a map $B^{-1}=\Delta^n\to F(U)$.
Lemma is proven.
\end{proof}

\subsubsection{}
Let $J$ be the collection of generating acyclic cofibrations. 
Let, furthermore, $\ol{J}$ denote the collection of maps which can be 
obtained as a countable direct composition of pushouts of coproducts of maps 
in $J$. Fibrations are defined as the maps satisfying RLP with respect to $J$. 

Repeating the reasoning of~\ref{end1:gen-com}, we prove this defines a model
structure on $\Dop X^{\pre}$. As in~\ref{end1:gen-com}, acyclic 
cofibrations in this model structure are retracts of elements of $\ol{J}$.

\Thm{simpl} is proven.

Note the following description of fibrations.

\subsubsection{}
\begin{prop}{s:fibrations}
A map $f:M\to N\in\Dop(X^{\pre})$ is a fibration if and only if the following
conditions are satisfied.
\begin{itemize}
\item{$f(U):M(U)\to N(U)$ is a Kan fibration for each $U\in X$.}
\item{For each hypercover $\epsilon:V_{\bullet}\to U$ the comutative
diagram
\begin{center}

\begin{diagram}
M(U)& \rTo & \Cech(V_{\bullet},M):=\Tot M(V_{\bullet})\\
\dTo&      &  \dTo                                    \\
N(U)& \rTo & \Cech(V_{\bullet},N)
\end{diagram} 
\end{center}
is homotopy cartesian.  }
\end{itemize} 
\end{prop}
\begin{proof}
A map $f:X\to Y$ of Kan simplicial sets is a weak equivalence iff
the space $X\times_YP(Y)\times_Y*$ is acyclic Kan for each point 
$*\to Y$. Here $P(Y)$ is the path space of $Y$. The rest of the proof
is a direct calculation.
\end{proof}

\end{document}